\documentclass[a4paper,10pt]{article}
\usepackage[plainpages=false]{hyperref}
\usepackage{amsfonts,latexsym,rawfonts,amsmath,amssymb,amsthm, mathrsfs, lscape}
\usepackage{verbatim}

\usepackage[all]{xy}
\usepackage{graphicx,psfrag}

\usepackage{array, tabularx}

\usepackage{setspace}
\usepackage[english]{babel}
\usepackage{color}

\usepackage{authblk}
\usepackage{enumerate}

\newtheorem{thm}{Theorem}[section]
\newtheorem{cor}[thm]{Corollary}
\newtheorem{lem}[thm]{Lemma}

\newtheorem{prop}[thm]{Proposition}

\newtheorem{conj}[thm]{Conjecture}

\theoremstyle{remark}
\newtheorem{rmk}[thm]{Remark}

\theoremstyle{definition}

\numberwithin{equation}{section}

\def \C {\textbf C}
\def \Z {\textbf Z}
\def \bR {\textbf R}
\def \R {\mathcal R}

\def \F {\mathcal F}
\def \I {\mathcal I}
\def \U {\mathcal U}
\def \T {T}

\def \ct {\mathfrak t}
\def \Q {\textbf Q}
\def \P {\mathbb P}
\def \H {\mathcal H}
\def \O {\mathcal O}
\def \M {\mathcal M}
\def \L {\mathcal L}
\def \m {\mathfrak m}
\def \t {\mathfrak t}
\def \K {\mathcal K}

\def \I {\text{I}}
\def \II {\text{II}}
\def \III {\text{III}}

\def \p {\partial}
\def \bp {\bar{\partial}}
\def \Vol {\text{Vol}}
\def \Spec {\text{Spec}}
\def \Hilb {\textbf{Hilb}}
\def \Aut {\text{Aut}}

\begin{document}

\title{Gromov-Hausdorff limits of K\"ahler manifolds and algebraic geometry, II}
\author{Simon Donaldson and Song Sun \thanks{The second author is partially supported by NSF grant DMS-1405832  and Alfred P. Sloan fellowship.}}

\date{}
\maketitle

\section{Introduction}
In this paper we study Gromov-Hausdorff limits of K\"ahler manifolds, in particular their singularities, and the connections with algebraic geometry. This is a continuation of our previous work \cite{DS}. 

For $\kappa>0$, let $\K(n, \kappa)$ be the class of $n$ dimensional compact polarized K\"ahler manifolds $(X, L, \omega, p)$, where $L$ is a Hermitian holomorphic line bundle over $X$ with curvature $-i\omega$ and $p$ is a chosen base point, which satisfy
 \begin{enumerate}[(1)]
 \item Einstein condition: 
 \begin{equation}\label{eqn1-1}
 Ric(\omega)=\lambda \omega
 \end{equation}
  with $|\lambda|\leq 1$;
 \item Local non-collapsing condition: for all $r\in (0, 1]$
 \begin{equation}\label{eqn1-2}
 \Vol(B(p, r))\geq \kappa r^{2n}.
 \end{equation}
  \end{enumerate}
  Notice these conditions are preserved if we rescale the metric by a factor that is greater than one. Given a sequence $(X_i, L_i, \omega_i, p_i)$ in $\K(n, \kappa)$,  from general theory by passing to a subsequence we obtain a (pointed) Gromov-Hausdorff limit $(Z, p)$, which is a length space.  By the regularity theory of Cheeger-Colding-Tian \cite{CCT}, we have a decomposition $Z=\mathcal R\sqcup \Sigma$, where $\mathcal R$ is an open connected smooth manifold endowed with a K\"ahler-Einstein metric,  and $\Sigma$ is a closed subset of $Z$ with real Hausdorff dimension at most $2n-4$.  Let $\O_\mathcal R$ be the structure sheaf of the  complex manifold underlying $\mathcal R$, and let $\iota: \mathcal R\rightarrow Z$ be the obvious inclusion map, then we define a sheaf on $Z$ by  $\O_Z=\iota_* \O_\mathcal R$. We have

\begin{thm} \label{thm1-1}
$(Z, \O_Z)$ is a normal complex analytic space. 
\end{thm}

The precise meaning of this statement, as well as further properties of $Z$,  will be discussed in detail in Section 2.    

\

 Although we expect the results below to hold in greater generality,  in this paper we will focus on the situation that is most closely related to our previous work \cite{DS}. For $V>0$  we denote by $\K_1(n, \kappa, V)$ the subset of $\K(n, \kappa)$ consisting of elements that satisfy the stronger assumptions  
 
\begin{enumerate}[(A)]
\item Normalized Einstein condition: 
$$Ric(\omega)=\lambda \omega$$
for $\lambda\in \{1, 0, -1\}$; if $\lambda=0$, we further assume $K_{X}$ is holomorphically trivial; 
\item Uniform non-collapsing condition:  
\begin{equation}\label{eqn1-3}
\Vol(B(q, r))\geq \kappa r^{2n}
\end{equation}
 for all $q\in X$ and $r\in(0, 1]$.
 \item Uniform volume bound: 
 \begin{equation}
 \Vol(X,\omega)\leq V.
 \end{equation}
\end{enumerate}

By the Bishop-Gromov volume comparison theorem,  (B) and (C) together are equivalent to a uniform diameter bound on $X$, and the latter is indeed a consequence of the Einstein condition when $\lambda=1$.    It is proved in \cite{DS} that the (polarized) Gromov-Hausdorff limit of a sequence of spaces in $\K_1(n, \kappa, V)$ is naturally  a normal projective variety. Theorem \ref{thm1-1} is an extension of this result.

Our main interest in this paper is on rescaled limits. For this purpose we let $\K(n, \kappa, V)$ be the set of polarized K\"ahler manifolds of the form $(X, L^{a}, a \omega, p)$ for some $(X, L, \omega, p)\in \K_1(n, \kappa, V)$ and $a\geq 1$.  Clearly  $\K(n, \kappa, V)$ is a subset of $\K(n, \kappa)$ so Theorem \ref{thm1-1} applies to Gromov-Hausdorff limits of spaces in $\K(n, \kappa, V)$.
  Let $(Z, p)$ be such a Gromov-Hausdorff limit.  We consider the family of spaces given by rescaling $(Z, p)$ by a factor $\sqrt{a}$ for a positive integer $a$.  Let $a\rightarrow\infty$, by passing to a subsequence we obtain limit spaces, called the  \emph{tangent cones} at $p$.  These can themselves be viewed as Gromov-Hausdorff limits of elements in $\K(n, \kappa, V)$, so by Theorem \ref{thm1-1} they are naturally complex analytic spaces. A fundamental result of Cheeger-Colding says that any tangent cone in this setting is also a metric cone, so is of the form $C(Y)$ for some compact metric space $Y$ (called the \emph{cross section}).  Let $R(C(Y))$ denote the ring of holomorphic functions on $C(Y)$ with polynomial growth at infinity. Then we have

\begin{thm} \label{thm1-2}
$R(C(Y))$ is finitely generated. Moreover,  $\Spec R(C(Y))$ is an affine algebraic variety that is complex analytically isomorphic to $(C(Y), \O_{C(Y)})$. 
\end{thm}

The proof of this theorem will be given in Section 2.3. We will also describe the cone structure on $C(Y)$, in terms of a ``grading" on  the ring $R(C(Y))$.

\begin{thm} \label{thm1-3}
There is a unique tangent cone at $p$. 
\end{thm}

We will prove this in Section 3. This means that any two tangent cones are isomorphic both as metric cones and as affine algebraic varieties, see also Remark \ref{rmk3-17} for more precise statement. 
 For general limits of Einstein manifolds, the uniqueness of tangent cones at a singular point is not known. In a recent work \cite{CM}, using a Lojasiewicz-Simon type argument, Colding-Minicozzi  proved the uniqueness under the  extra assumption that there is one tangent cone with smooth cross section. Our approach  is very different from this in that we exploit the complex geometry in a crucial way and the above theorem does not require the smoothness of the cross section. In Section 3 we also make some progress towards an algebro-geometric description of the tangent cone. In particular, we will study the relation between the grading on $R(C(Y))$ and the filtration on the local ring of germs of holomorphic functions at $p$ defined by the limit metric. 
 
 \
 
When the above limit space $(Z, p)$ is non-compact, i.e. when the rescaling factors $a_i$ tend to infinity, we can ask about the algebraicity of $Z$. Let $R(Z)$ be the ring of holomorphic functions on $Z$ with polynomial growth at infinity.  Then we have

\begin{thm} \label{thm1-4}
$R(Z)$ is finitely generated. Moreover,  $\Spec(R(Z))$ is an affine algebraic variety that is complex analytically isomorphic to $(Z, \mathcal O_Z)$. 
\end{thm}

This is an extension of Theorem \ref{thm1-2}. The proof is given in Section 3.4. It involves the study of tangent cones at infinity, for which we will obtain results analogous to Theorem \ref{thm1-2} and \ref{thm1-3}. By our assumption $Z$ is also endowed with a Ricci-flat K\"ahler metric. When $Z$ is smooth, it is an \emph{asymptotically conical Calabi-Yau manifold},  which has been well-studied  recently (see for example \cite{CH}).  Theorem \ref{thm1-4} can also be compared with \cite{Liu}, where a similar result is proved for complete K\"ahler manifolds with non-negative bisectional curvature and maximal volume growth. 

\

In the appendix we will prove an extension of the Futaki and Matsushima theorem to singular Ricci-flat K\"ahler cones, which is used in the proof of Theorem \ref{thm1-3} and \ref{thm1-4}.  Our arguments follow the corresponding proof for $\Q$-Fano varieties in \cite{CDS3}.

\

The main application of our results in this paper is to the study of K\"ahler-Einstein metrics with positive Ricci curvature (i.e. the \emph{Fano} case), in which case the non-collapsing condition holds automatically.  For K\"ahler-Einstein metrics with negative or zero Ricci curvature (i.e. the \emph{General Type} or \emph{Calabi-Yau} case, respectively),   it is an interesting question to understand the algebro-geometric meaning of the non-collapsing condition. There are recent results along this direction, see for example  \cite{RZ, Tosatti, Song}.

\

\textbf{Acknowledgements:} We are grateful to Mark Haskins, Weiyong He, Hans-Joachim Hein, Robert Lazarsfeld and Jason Starr for helpful discussions related to this work.

 \section{Complex structure on  Gromov-Hausdorff limits}

\subsection{Proof of Theorem 1.1}

We first recall the notion of polarized Gromov-Hausdorff convergence introduced in \cite{DS}. Fix $n$ and $\kappa>0$, suppose we are given a sequence of objects $(X_i, L_i,  \omega_i, p_i)$ in $\K(n, \kappa)$.  Then by passing to a subsequence we obtain a \emph{polarized limit space} $(Z, p, g_\infty, J_\infty, L_\infty, A_\infty)$, which consists of the Gromov-Hausdorff limit metric space $(Z, p)$, together with a smooth Riemannian metric $g_\infty$ and a  compatible complex structure $J_{\infty}$ on the regular set $\R$ with  K\"ahler form $\omega_{\infty}$, a Hermitian line bundle $L_{\infty}$ over $\R$, and a smooth connection $A_{\infty}$ on $L_{\infty}$ whose curvature is $-i\omega_{\infty}$ (The difference from the definition in \cite{DS} is that here we assume the metrics satisfy the Einstein equation so the limiting geometric structures are all smooth over $\R$).

The meaning of the convergence is as follows. For any $R>0$, we can fix a metric $d_i$ on the disjoint union $B(p_i, R)\sqcup B(p, R)$ such that $B(p_i, R)$ and $B(p, R)$ are both $\epsilon_i$-dense, and $d_i(p_i, p)\leq \epsilon_i$ with $\epsilon_i\rightarrow 0$. Moreover for any $\delta>0$, and any compact subset $K\subset  B(p, R)\cap \R$  we can find for large enough $i$ open embeddings $\chi_{i}$ of an open neighbourhood of $K$ into $B(p_i, R)$, and bundle isomorphisms $\hat\chi_i: L_\infty\rightarrow \chi_i^*L_i$,  such that $d_i(x, \chi_{i}(x))\leq \delta$ for all $x\in K$, and $(\chi_i^*g_{i}, \chi_i^*J_i, \chi_i^*A_i)$ converges smoothly over $K$ to $(g_{\infty}, J_\infty, A_\infty)$. Here $J_i$ is the complex structure on $X_i$ and $A_i$ is the Chern  connection on $L_i$ with curvature $-i\omega_i$.

\

Now let $(Z, p)$ be a limit space. Fix $R>1$, and fix a metric $d_i$ on $B(p_i, R)\sqcup B(p, R)$ which realizes the polarized Gromov-Hausdorff convergence. Let $\mathcal O$ be the sheaf of rings on $B(p, R)$  induced by the presheaf on $B(p, R)$, which assigns each $\Omega\subset B(p, R)$ the ring of functions on $\Omega$ that are limits of holomorphic functions over certain domains in $X_i$, in the obvious sense. From the definition $\O$ depends on $R$ and the choice of $d_i$, but eventually we will prove that $\O$ agrees with the restriction of $\O_Z$ defined in the introduction, so it in fact does not depend on any choices. 

The overall idea to prove Theorem \ref{thm1-1} is similar to the one we used in the proof of Theorem 2 in \cite{DS}. 
 As discussed in \cite{DS}, \cite{CDS2} all tangent cones of $Z$ are ``good" so that we can apply the H\"ormander technique to construct holomorphic sections.   Recall in \cite{CDS2}  we have achieved the following. 
 
 \begin{prop} \label{prop2-1}
 There are $k$, $C$, $N$, $l_1, \cdots, l_N$, and $\rho_1, \rho_2\in (0, 1)$ with $\rho_1>\rho_2$,  such that the following holds
 
\begin{enumerate}
\item[(1)] For $i$ sufficiently large there is  a holomorphic section $s_i$ of $L_i^{k}\rightarrow X_i$ such that  $|s_i(x)|\geq 1/2$ when $d_i(p, x)\leq \rho_1$, and $||s_i||_{L^{2}}\leq (2\pi)^n+1$, where the norm is measured with respect to the metric $k\omega_i$;
\item[(2)] For $j=1, \cdots, N$, a holomorphic section $\sigma_{j, i}$ of $L_i^{kl_j}\rightarrow X_i$ for some integer $l_j\geq 1$;
\item[(3)] The corresponding map $F_i: D_i\rightarrow \C^N$,  with the $j$-th component  given by $\sigma_{j, i}/s_i^{l_j}$, satisfies $|F_i(x)|_*> 1/2$ when $d_i(p, x)=\rho_1$, and $|F_i(x)|_*\leq 1/100$ when $d_i(p, x)\leq \rho_2$. Here $|\cdot|_*$ is the sup  norm on $\C^N$, and $D_i$ is a domain in $X_i$ containing all the points with $d_i(p, x)\leq \rho_1$;
\item[(4)] $|\nabla F_i|\leq C. $
\end{enumerate}
\end{prop}
 
 Let $B$ be a Euclidean ball in $\C^N$ that is contained in the ball of radius $1/4$ in the $|\cdot|_*$ norm, and let $\Omega_i$ be the pre-image of $B$ under $F_i$. Item (3) implies that for any $x\in B$,  the fiber $F_i^{-1}(x)$ is a compact analytic set and by item (1) the ample line bundle $L_i^{ka_i}$ is trivial over $F_i^{-1}(B)$, therefore $F_i$ is a finite map from $\Omega_i$ onto an analytic set $W_i$ in $B$. Item (4) means that the volume of $W_i$ measured by the induced metric from $\C^N$ is uniformly controlled by the volume of a ball of radius $\rho_1$ in $X_i$ (with respect to the metric $g_i$), and the latter is uniformly bounded by the Bishop-Gromov volume comparison theorem. So by passing to a subsequence we may assume $W_i$ converges to a limit $W$, which is an analytic set in $B$. We endow $W$ with the reduced analytic structure. (4) also implies that we can take the limit of $F_i$ and obtain a Lipschitz map $F$ from an open neighborhood $\Omega$ of $p$ onto $W$. This induces an injective sheaf map $F^*: \mathcal O_W\rightarrow F_* \O$, where $\mathcal O_W$ is the sheaf of holomorphic functions on $W$.

 \begin{prop} \label{prop2-2}
 The following can be achieved:
 
(A). For any $q_1, q_2\in \Omega$, there are an integer $r$ and $\epsilon>0$,  and for $i$ large there are holomorphic sections $\tau_{1, i}, \tau_{2, i}$ of $L^{kr}\rightarrow X_i$ so that the functions $\tau_{1, i}/s_i^r, \tau_{2, i}/s_i^r$ converge to functions in  $\mathcal O(\Omega)$ which separate $B(q_1, \epsilon)$ and $B(q_2, \epsilon)$;
 
 (B).  For any point $q\in \Omega\cap\R$, there is an integer $r$, and for $i$ large there are $n$ holomorphic sections $(\tau_{1,i}, \cdots, \tau_{n, i})$ of $L^{kr}\rightarrow X_i$ so that $(\tau_{1, i}/s_i^{r}, \cdots, \tau_{n, i}/s_i^{r})$ converge to functions in $\mathcal O(\Omega)$ that define an embedding of an open neighborhood of $q$ into $\C^n$.   
 
 \end{prop}
 
The proof of these is exactly the same as that of Proposition 4.6 and Proposition 4.7 in \cite{DS}, by constructing Gaussian holomorphic sections around two different points separately, and by constructing holomorphic sections which vanish at one point but with non-vanishing derivative along any prescribed tangent direction.

Given a function $f\in \mathcal O(\Omega)$,  we could add it as a new component and obtain a map $F'=(F, f): \Omega\rightarrow \C^{N+1}$. By definition, $f$ is the limit of holomorphic functions $f_i$ defined over some open subset in $X_i$. By the gradient estimate for holomorphic functions (see for example Proposition 2.1 in \cite{DS}), $|\nabla f_i|$ is locally uniformly bounded. This implies that the image $W'$ of $F'$ is a local complex analytic set in $\C^{N+1}$. Moreover the projection map $\pi:W'\rightarrow W$ is  finite. By (A) and (B) we may add finitely many  components so that the map $F'$ is one-to-one from some open subset $D$ in $\Omega\cap \R$ onto an open subset of the smooth part of $W$, and the pre-image of $F'(D)$ is exactly $D$. Without loss of generality  we may assume $F$ already meets this property. 

\begin{prop} \label{prop2-3}
By adding finitely many functions from $\O(\Omega)$ and by slightly shrinking $\Omega$, we may assume $F$ is a homeomorphism and maps $\Omega\cap\R$ into the smooth part of $W$. 
\end{prop}

If we add another function in $\O(\Omega)$ as a new component, the projection map will be generically one-to-one, so in particular we have the induced inclusion of sheaves of rings $\O_W\hookrightarrow \pi_*\O_{W'}\hookrightarrow \M_W$, where $\M_W$ is the sheaf of meromorphic functions on $W$. Indeed,  $\pi_*\O_{W'}$ is  a coherent subsheaf of $\widehat{\O_{W}}$, the normalization of $\O_W$.  By general theory of complex spaces (c.f. \cite{Remmert}, Section 11.5) we have a Noether property, that is, by adding finitely many functions from $\O(\Omega)$ and by slightly shrinking $\Omega$,  we may eventually achieve a maximal subsheaf, say $\pi_*\O_{W'}\subset\widehat{\O_W}$.   Again,  without loss of generality we may assume $\O_{W}$ is already maximal, then we have $\O(W)=F_*(\O(\Omega))$.  Since the functions constructed from (A) and (B) clearly lie in $\O(\Omega)$,  we see  that the map $F$ is a homeomorphism onto $W$ and it maps $\Omega\cap\R$ into the smooth part of $W$.

\

Now we may run the same arguments locally.  Using the fact that $\O_{W,p}$ is a Noetherian ring, by adding functions in $\O_p$ and by shrinking $\Omega$ if necessary, we may assume that  $\O_{W,p}= F_*(\O_p)$. 

\begin{prop}\label{prop2-4}
By further shrinking $\Omega$ if necessary we may assume $W$ is normal, and the map $F^*: \O_W\rightarrow F_*(\O|_\Omega)$ is bijective.
\end{prop}

By the openness of normal locus (\cite{Remmert}, Theorem 14.4) it suffices to show $W$ is normal at $F(p)$. This is a local property, so without loss of generality we may assume $W_i$ and $W$ are analytic subsets of a Euclidean ball $B$ in $\C^N$, and we need to prove that a bounded holomorphic function $f$ defined over the smooth part of $W$ extends to a holomorphic function over a neighborhood of $F(p)$ in $W$. By the above discussion, it suffices to prove that any bounded holomorphic function $f$ defined over $\overline{\Omega}\cap\R$  extends to a function in $\O(\Omega)$. For this purpose we need to use a local version of  the H\"ormander $L^2$ estimate.  The following lemma is well-known, see for example \cite{Demailly}, Theorem 6.1. 
 
\begin{lem}\label{lem2-5}
Let $Y$ be a complex manifold which admits a complete K\"ahler metric. Let $\omega$ be an arbitrary K\"ahler form on $Y$,   and $L$ be a holomorphic line bundle over $Y$ endowed with a Hermitian metric $h$ whose curvature satisfies $i\Theta_h \geq c\omega$ for some $c>0$. Let $f$ be an $L$-valued $(n, q)$ $(q\geq1)$ form with $\bp f=0$, then there exists an $L$-valued $(n, q-1)$ form $u$ with  $\bp u=f$, and 
$$||u||_{L^2}^2\leq (cq)^{-1}||f||_{L^2}^2. $$
\end{lem}
\

 Notice  each $W_i$ is an analytic set in $B$, so is Stein. Since the map $F_i$ is finite,  it is easy to see that $\Omega_i$ admits  a complete K\"ahler metric. We also choose a big number $r$ so that on $\Omega_i$ the curvature of the  line bundle $L_i^{kr}\otimes K_{X_i}^{-1}$ is bigger than $\omega_i$, where the metric on $K_{X_i}^{-1}$ is defined by $\omega_i^n$. Now as in \cite{DS} we fix a sequence $\eta_i\rightarrow0$. Using the fact that the singular set $\Sigma$ has Hausdorff dimension strictly less than $2n-2$, we can choose a sequence of good cut-off functions $\beta_i$ on $\Omega$ so that $\beta_i$ is supported  in the complement of a neighborhood of $\Sigma\cap \Omega$,  $\beta_i=1$ outside the $\eta_i$-neighborhood of $\Sigma$ and $||\nabla \beta_i||_{L^2}\leq \eta_i$. Given a non-zero bounded holomorphic function $f$ defined over $\overline{\Omega}\cap \mathcal R$, 
we can use the maps $\chi_i$ to graft $\beta_i  f$ into $\Omega_i$ and obtain a smooth section $\sigma_i=(\chi_i^{-1})^*(\beta_i f) s_i^{\otimes r}$ of  $L_i^{kr a_i}$ over $\Omega_i$, with $||\bp \sigma_i||_{L^2}\rightarrow 0$. We may view $\sigma_i$ naturally as a $L_i^{kra_i}\otimes K_{X_i}^{-1}$-valued $(n, 0)$ form.  By Lemma  \ref{lem2-5}
we can solve $\bp \tau_i= \bp \sigma_i$ with $||\tau_i||_{L^2}\leq ||\bp\sigma_i||_{L^2}$.

Let $f_i=(\sigma_i-\tau_i)/s_i^{\otimes r}$, then $f_i$ is a holomorphic function on $\Omega_i$, and as in \cite{DS}  we obtain a uniform $L^\infty$ estimate on $f_i$ and $|\nabla f_i|$, where the constants depend only on the distance to the boundary of $\Omega_i$. So we can take a limit $f_\infty\in \O(\Omega)$ by passing to a subsequence. On a ball in $\Omega_i$ that has a fixed distance away from $\Sigma$, we then obtain a uniform estimate on $|\nabla \tau_i|$, and this together with the fact that $||\tau_i||_{L^2}\leq\epsilon_i$ implies that $|\tau_i|$ tends to zero uniformly on any compact subset of $\Omega\cap \mathcal R$.  Therefore  $f_\infty=f$ on the whole $\Omega\cap\R$, hence can be viewed as an extension of $f$ to $\Omega$.  

 Now without loss of generality we may assume $W$ itself is normal. 
We need to show $F^*: \O_W\rightarrow F_*(\O|_\Omega)$ is bijective. It suffices to prove the surjectivity. Given any $q\in \Omega$, a holomorphic function $f$ defined on a neighborhood $U$ of $q$  determines a holomorphic map $F': U\rightarrow\C^{N+1}$, and the projection map $\pi: F'(U) \rightarrow F(U)$ is generically one-to-one, and $f$ becomes holomorphic on $F'(U)$. On the other hand,  the normality of $W$ implies that $\pi$ is a holomorphic equivalence, so $f$ is holomorphic on $F(U)$.  This completes the proof of Proposition \ref{prop2-4}.\\

Notice by normality $\O_Z|_{\Omega}=(\iota_{*}\O_\R)|_{\Omega\cap \R}=\O|_{\Omega}$. Therefore $F$ also induces an isomorphism between the ringed spaces $(\Omega, \O_Z|_\Omega)$ and $(W, \O_W)$, and the same holds in a neighborhood of  any point in $Z$. This then endows $(Z, \mathcal O_Z)$ with the structure of a normal complex space, in the usual sense, and thus finishes the proof of  Theorem 1.1. \\

\begin{rmk}
 The above arguments make use of some general language of complex analytic spaces, and are essentially equivalent to the  approach used in \cite{DS}. 
  \end{rmk}

\subsection{Further results}

We first clarify the precise notion of Gromov-Hausdorff topology we shall use in the remaining part of this paper. We go back to the setting at the beginning of Section 2.1, where we discuss a polarized limit space $(Z, p, g_\infty, J_\infty, L_\infty, A_\infty)$. Since the Hermitian line bundle and the connection enter our discussion only when we apply the H\"ormander construction in \cite{DS} (more specifically the construction of holomorphic sections in Proposition \ref{prop2-1} and \ref{prop2-2}),  and since they are not the geometric objects that we are interested in later, we will mostly ignore them. From now on, we will simply call $(Z, p, g_\infty, J_\infty)$, or $(Z, p)$ when there is no confusion caused, a \emph{Gromov-Hausdorff limit}. By abusing notation we will also denote by $\K(n, \kappa)$ the class of the underlying (non-polarized) K\"ahler manifolds of elements in $\K(n, \kappa)$ defined in the introduction. Let $\overline{\K(n, \kappa)}$ be the class of all Gromov-Hausdorff limits of elements in $\K(n, \kappa)$.  It is understood that an element $(Z, p)\in \overline{\K(n, \kappa)}$ is always endowed with some limit polarization, but is in general not unique.

 The discussion in Section 2.1 defines a notion of convergence in $\overline{\K(n, \kappa)}$ (by forgetting about the line bundle and connection) which, by  general construction,  yields topology on  $\overline{\K(n, \kappa)}$. This refines the standard Gromov-Hausdorff topology on metric spaces, and this is what we mean by  \emph{Gromov-Hausdorff topology} in the rest of this paper. 

A basis of this topology can be constructed as follows. Given a positive integer $j$  and $(Z, p)\in \overline{\K(n, \kappa)}$, we define a neighborhood  $N_j(Z, p)$ to be the set of all spaces $(Z', p')\in \overline{\K(n, \kappa)}$ which satisfy the following properties

\begin{itemize}
\item There is a metric $d$ on $\overline{B(p, j)}\sqcup \overline{B(p', j)}$, such that $d(p', p)< j^{-1}$ and $\overline{B(p, j)}$ and $\overline{B(p', j)}$ are both $\epsilon$-dense for some $\epsilon<j^{-1}$;
\item Denote by $U_j$ the complement of the $j^{-1}$ neighborhood of the singular set in $B(p, j)$. Then there is a smooth embedding $\chi$ of an open neighborhood of $U_j$ into the smooth part of  $B(p', j)$, such that $d(x, \chi(x))< j^{-1}$ for all $x\in U_j$, and
 $$||\chi^*g'-g||_{C^j(U_j)}+||\chi^*J'-J||_{C^j(U_j)}< j^{-1},$$
 where the norm is computed with respect to the metric $g$.  
\end{itemize}
Then the collection of the neighborhoods $N_j(Z, p)$ for all $(Z, p)$ and all integers $j$ is a basis of the Gromov-Hausdorff topology on $\overline{\K(n, \kappa)}$.

\begin{lem} \label{lem2-0}
The Gromov-Hausdorff topology on $\overline{\K(n, \kappa)}$ is compact, Hausdorff, and has a countable basis. 
\end{lem}

This is not difficult to prove but since we can not find a reference for the precise statement in the literature we give a proof here.   By a contradiction argument, it is easy to see that for each fixed $j$, there are at most finitely many disjoint neighborhoods of the form $N_j(Z, p)$. Let $N_j(Z_{j,\alpha}, p_{j, \alpha})$ ($\alpha=1, \cdots, c(j)$) be a maximal disjoint set of such neighborhoods. Now we claim the countable family of open subsets $\{N_j(Z_{k, \alpha})|\alpha=1, \cdots, c(k), j\in \Z_{>0}\}$ form a basis of the topology. To see this, given an open set $N$ in $\overline{\K(n, \kappa)}$ and a point $(Z, p)\in N$, by our above choice, for each $j$, we can find $\alpha_j$ such that $N_j(Z, p)$ and $N_j(Z_{j, \alpha_j}, p_{j, \alpha_j})$ has non-empty intersection. Then it is easy to see that for $j$ sufficiently large we have $N_{100j}(Z, p)\subset N_{10j}(Z_{j, \alpha_j}, p_{j, \alpha_j})\subset N_{j}(Z, p)\subset N$. This proves the claim, and hence the topology has a countable basis. 

The Hausdorff property is clear. It remains to prove the compactness. Suppose otherwise, we may find an open cover $\mathcal U$ of $\overline \K(n, \kappa)$ which does not admit any finite sub-cover. Using the countable basis constructed above we may choose a countable sub-cover, say $\mathcal V=\{V_1, V_2, \cdots\}$. By assumption, for each $k$, there is a $(Z_k, p_k)$ which does not belong to $V_j$ for any $j\leq k$.  Since each $(Z_k, p_k)$ is the Gromov-Hausdorff limit of a sequence of spaces in $\K(n, \kappa)$, by a diagonal sequence argument,  we may pass to a subsequence and assume $(Z_k, p_k)$ converges to a limit $(Z, p)$. Now since $\mathcal V$ is a cover,  $(Z, p)\in V_{k_0}$ for some $k_0$. It follows that  for $k$ sufficiently large $(Z_k, p_k)\in V_{k_0}$. Contradiction.  
 
\

For our later purposes we need to extend the discussion of Section 2.1 uniformly to  $\overline{\K(n, \kappa)}$. 

\begin{lem}\label{lem2-7}
There are  $\lambda_1, \lambda_2\in (0, 1)$ with $\lambda_1>\lambda_2$ and $C>0$ depending only on $n$ and $\kappa$ such that given $(Z, p)\in \overline{\K(n, \kappa)}$,  there is an open set $D$ in $Z$ that contains the closure of the ball $B_{\lambda_1}(p)$, and a holomorphic map $F$ from $D$ to $\C^N$ such that 
\begin{itemize}
\item $|F(x)|_*\geq1/2$ when $d(x, p)=\lambda_1$;
\item  $|F(x)|_*\leq 1/100$ when $d(x, p)\leq \lambda_2$;
\item $|\nabla F(x)|\leq C$ for all $x$ in $B_{\lambda_1}(p)$. 
\end{itemize}
\end{lem}

This follows directly from Proposition \ref{prop2-1} and a contradiction argument. 

\begin{prop} \label{prop2-8}
 Suppose $f$ is holomorphic function defined on $B_{\lambda_1}(p)$. Then there is a neighborhood $\mathcal U$ of $(Z, p)$ in $\overline{\K(n, \kappa)}$, such that for any $(Z', p')\in \overline{\K(n, \kappa)}$, there is a holomorphic function $f'$ defined on  $B_{\lambda_2}(p')$, such that  $f'$ converges to  $f$ uniformly over $B_{\lambda_2}(p)$ as $(Z', p')$ converges to $(Z, p)$.
\end{prop}

\begin{rmk} \label{rmk2-10}
Notice the precise notion of convergence of holomorphic functions in our context depends on the choice of metric on the disjoint union $Z'\sqcup Z$ realizing the  Gromov-Hausdorff convergence. In general the limit will be only well-defined up to an isomorphism of $B_{\lambda_2}(p)$ (a holomorphic isometry) that fixes $p$. So the precise convergence should be understood modulo such an isomorphism. 
\end{rmk}

By Lemma \ref{lem2-7}, we can apply Lemma \ref{lem2-5} and the discussion following it to find the neighborhood $\mathcal U$, and a holomorphic function $f'$ on $B_{\lambda_2}(p')$ for $(Z', p')\in \mathcal U \cap \K(n, \kappa)$, with a uniform $L^\infty$ estimate on $|f'|$ and $|\nabla f'|$. Then  by taking limits of these functions we also find correspondingly the holomorphic functions for all $(Z', p')\in \mathcal U$.

\

\begin{prop} \label{prop2-10}
There are universal constants $K_0$, $K_1$ depending only on $n$, $\kappa$ and $r$, so that for any holomorphic function $f$ defined on a ball $B$ of radius $r$ around $p$ in a limit space $Z$, we have  
$$|f(p)|\leq K_0 |f|_{L^2(B)}; $$
$$|\nabla f(p)|\leq K_1|f|_{L^2(B)}.$$
Here the second estimate is understood in the Lipschitz sense. 
\end{prop}

A general way to prove this is to adapt the usual Moser iteration technique directly to the possibly singular space $Z$. In our  case, we can apply Proposition \ref{prop2-8} and the fact that the estimate is well-known in the case when $Z$ is in $\K(n,\kappa)$ (c.f. Proposition 2.1 in \cite{DS}).

\begin{prop}  \label{prop2-11}
Let  $F: B_{\lambda_1}(p)\rightarrow \C^N$ be a holomorphic embedding. Then we may find $\mathcal V\subset \mathcal U$, and for any $(Z', p')\in \mathcal V$, a holomorphic map $F': B_{\lambda_2}(p')\rightarrow\C^N$, that is generically one-to-one  (in particular $F'$ is a normalization map onto its image), and as $(Z', p')$ converges to $(Z, p)$,  the image  $F'(B_{\lambda_2}(p'))$ converges to $F(B_{\lambda_2}(p))$ as local complex analytic sets in $\C^N$. 
\end{prop}

The construction of $F'$ follows from Proposition \ref{prop2-8}.  By compactness of $\overline{\K(n, \kappa)}$,  to prove that $F'$ is generically one-to-one,  it suffices to show that if $(Z_i, p_i)$ converges to $(Z, p)$ then $F_i$ is generically one-to-one for sufficiently large $i$. Fix a metric of $B_{\lambda_1}(p_i)\sqcup B_{\lambda_1}(p)$ that realizes the Gromov-Hausdorff convergence. Choose a ball $B$ with closure contained in the regular part of $B_{\lambda_2}(p)$. Then we may find corresponding balls $B_i$ in $B_{\lambda_2}(p_i)$ that converge to $B$. By Colding's volume convergence theorem and Anderson's volume gap theorem it follows that for $i$ large $B_i$ is contained in the regular part of $B_{\lambda_2}(p_i)$. From the definition of Gromov-Hausdorff convergence, by varying $B_i$ slightly we may identify $B_i$ with $B$ using a diffeomorphsim $\chi_i$, under which $F_i$ converges smoothly to $F$. Hence $F_i$ is an embedding on $B_i$. Now using the injectivity of $F$ it is easy to see that for $i$ large on the image of the half ball $\frac{1}{2}B_i$, $F_i$ is one-to-one, in particular, $F_i$ is generically one-to-one. 

\begin{rmk}
 Proposition \ref{prop2-11} will be sufficient for our purpose in this paper. In general we expect that for $i$ large $F_i$ is indeed a holomorphic embedding,  which will allow us to say that $Z_i$ converges to $Z$ locally as analytic subsets  in some  $\C^N$, or in other words, locally $Z_i$ is a \emph{deformation} of $Z$.  Comparing the results in the compact case  \cite{DS} that relate the Gromov-Hausdorff convergence to flat convergence in the Hilbert scheme, we also expect that in general the convergence is  \emph{flat} in a certain sense.  To our knowledge such a theory  has not yet been developed, and we leave this for future work.
\end{rmk}

\begin{prop} \label{prop2-13}
The metric singular set of $(Z, p)$ agrees with the complex analytic singular set.
\end{prop}

By Proposition \ref{prop2-3} it suffices to show that if $p$ is a smooth point in the complex analytic sense, then the limit metric is smooth in a neighborhood of $p$. Choose a holomorphic embedding of a ball $B$ around $p$ into $\C^n$. By a rescaling we may assume $B=B_{\lambda_1}(p)$. Then by Proposition \ref{prop2-11} for $i$ large enough we may find a holomorphic map $F_i: B_{\lambda_2}(p_i)\rightarrow \C^n$, that is generically one-to-one, hence is a holomorphic equivalence onto its image, and $F_i$ converges to $F$. By making $B$ even smaller we may view the K\"ahler-Einstein metric $\omega_i$ as a K\"ahler metric on a fixed Euclidean ball $\textbf{B}$ in $\C^n$. Moreover, we can write $\omega_i=i\p\bp\phi_i$, with $\phi_i=-k^{-1}\log|s_i|^2$, where $s_i$ is the holomorphic section constructed in Proposition \ref{prop2-1}.   As in \cite{CDS2},  we have $|\phi_i|\leq C$ and $\omega_i\geq C^{-1}\omega_{Euc}$ for some $C>0$. The K\"ahler-Einstein equation takes the form
$$\det(i\p\bp\phi_i)=e^{-\lambda \phi_i}|U_i|^2,  $$
where $U_i$ is a non-vanishing holomorphic function. The bound $\omega_i\geq C^{-1}\omega_{Euc}$ implies that $|U_i|^{-1}$ is uniformly bounded. Since the volume of $\textbf{B}$ with respect to $\omega_i$ is uniformly bounded we obtain a $L^2$ bound on $U_i$. So in the smaller ball, say $\frac{3}{4}\textbf{B}$, we know $|U_i|$ is also uniformly bounded.  This implies $\omega_i$ and $\omega_{Euc}$ are uniformly equivalent. Then we can apply the standard Evans-Krylov theory to conclude that $\phi_i$ has a uniform $C^{2,\alpha}$ bound on $\frac{1}{2}\textbf{B}$, and standard bootstrapping yields higher derivative bound. So $\omega_i$ converges to a smooth K\"ahler-Einstein metric $\omega_\infty$ in $\frac{1}{2}\textbf{B}$.

\begin{rmk}
As in \cite{DS}, the above argument also proves that there is a \emph{weak} K\"ahler-Einstein metric on $Z$ in the sense of pluri-potential theory, with continuous local potential.
\end{rmk}

\

We finish this subsection with a lemma on the convergence of holomorphic functions, that will be used later. 
 Suppose a sequence $(Z_i, p_i)\in \overline{\K(n, \kappa)}$ converges to $(Z, p)$. Suppose $B_i$ is a ball in $Z_i$ that converges to a ball $B$ in $Z$. Given a sequence of holomorphic functions $f_i$ on $B_i$ with $||f_i||_{L^2(B_i)}$ uniformly bounded, then by the estimate in Lemma \ref{prop2-10},  we know $f_i$ converges (by passing to a subsequence) to a holomorphic function $f$ on $B$,  and  the convergence is uniform over any compact subset of $B$. In this case we say \emph{$f_i$ converges weakly to $f$.  }
 
 From our definition of Gromov-Hausdorff convergence,  any domain $\Omega$ with $\bar\Omega\subset B\cap\mathcal R$ is the smooth limit of domains $\Omega_i$ in $B_i$, so  we always have 
 $$||f||_{L^2(B)}\leq \liminf_{i\rightarrow\infty} ||f_i||_{L^2(B_i)}.$$
We say \emph{$f_i$ converges \emph{strongly} to $f$} if 
 $$||f||_{L^2(B)}= \lim_{i\rightarrow\infty} ||f_i||_{L^2(B_i)}.$$

\begin{lem} \label{lem2-15}

(1). If $f_i$ converges uniformly to $f$, then $f_i$ converges strongly to $f$. \\
(2). Let $\lambda_1, \lambda_2$ be given as in Lemma \ref{lem2-7}.  Suppose $f_i$ converges strongly to $f$, and $g_i$ converges weakly to $g$. If $f$ extends to a holomorphic function over the ball $B'=\lambda_1\lambda_2^{-1}B$ (the same center but with radius multiplied by $\lambda_1\lambda_2^{-1}$), then $$\int_{B} f\bar g=\lim_{i\rightarrow\infty}\int_{B_i} f_i\bar g_i. $$
\end{lem}

Given a domain $\Omega\subset B\cap \mathcal R$, we can find $\Omega_i\subset B_i$ which converges smoothly to $\Omega$. Then we have $||f_i||_{L^2(\Omega_i)}$ converges to $||f|_{L^2(\Omega)}$. On the other hand, by Colding's volume convergence theorem,  we know $\Vol(B_i\setminus \Omega_i)$ converges to $\Vol(B\setminus\Omega)$, which can be made as small as we like, 
 since the singular set in $B_\infty$ has zero $n$-dimensional Hausdorff measure. This proves the first item.

To prove the second item we first use Proposition \ref{prop2-8} to find a holomorphic function $h_i$ on $B_i$ that converges uniformly to $f$. 
Then we claim
 \begin{equation} \label{eqn2-1}
 \lim_{i\rightarrow\infty} \int_{B_i} h_i \bar g_i=\int_{B} f \bar g.
 \end{equation}
   To see this, let $rB_i$ and $rB$ be the balls with the same center as $B_i$ and $B$ respectively, and with radius multiplied by $r$. For any fixed $r<1$, since $g_i$ converges uniformly to $g$ on $rB$,  by item (1),
$\int_{rB_i}  h_i \bar g_i$ converges to $\int_{rB} f \bar g. $
On the other hand, we have
$$|\int_{B_i\setminus  rB_i} h_i  \bar g_i|\leq |h_i|_{L^\infty}|g_i|_{L^2(B_i)} \sqrt{\Vol(B_i\setminus  rB_i)}$$
As $r\rightarrow 1$, the right hand side tends to $0$ uniformly for all $i$. This proves the claim.
Now we write
\begin{equation} \label{eqn2-2}
\int_{B_i} f_i \bar g_i=\int_{B_i} (f_i-h_i)\bar g_i+\int_{B_i} h_i \bar g_i.
\end{equation}
Notice
\begin{equation*}
 \int_{B_i} |f_i-h_i|^2 =\int_{B_i} |f_i|^2+|h_i|^2-2Re(h_i\bar f_i)
 \end{equation*}
By assumption the first term converges to $||f||_{L^2(B)}^2$. By item (1), the second term also converges to $||f||_{L^2(B)}^2$.   Applying (\ref{eqn2-1}) with $g_i$ replaced by $f_i$, we see the last term converges to $-2||f||_{L^2(B)}^2.$  These imply $||f_i-h_i||_{L^2(B_i)}$ converges to zero. Since $||g_i||_{L^2(B_i)}$ is uniformly bounded by assumption, it follows that the first term in (\ref{eqn2-2}) converges to zero. Therefore by (\ref{eqn2-1}) again we obtain the conclusion.

\subsection{Proof of Theorem \ref{thm1-2}}
 Let $(Z, p)$ be a Gromov-Hausdorff limit  of a sequence of spaces in $\K(n, \kappa, V)$.   Let $C(Y)$ be a tangent cone at $p$. Then we know $C(Y)$ is in $\overline{\K(n, \kappa)}$, so 
 Theorem \ref{thm1-1} already proves that $C(Y)$ has the structure of a normal complex analytic space.  The main new ingredient in Theorem \ref{thm1-2} is the algebraicity. We will make use of an idea due to Van Coevering \cite{VC}, who proved essentially the same result for a K\"ahler cone with smooth cross section.

Using the metric cone structure we may write the smooth part of $C(Y)$ as $C(Y^{reg})$, where $Y^{reg}\subset Y$ is a smooth $2n-1$ dimensional manifold. The K\"ahler-Einstein condition implies that $C(Y^{reg})$ is Ricci-flat K\"ahler and $Y^{reg}$ is Sasaki-Einstein with Ricci curvature $2n-2$. Let $\xi=J (r{\p r})$, where $r$ is the distance function to the vertex $O$, and $J$ is the complex structure on $C(Y^{reg})$. By a simple local calculation it is easy to see that $\xi$ is holomorphic and Killing on $C(Y^{reg})$. 

\begin{lem} \label{compact torus action}
$\xi$ generates a holomorphic isometric action of a compact torus $\mathbb T$ on $C(Y)$. 
\end{lem}

We choose a neighborhood $\Omega$ of $O$ and a holomorphic embedding $F: \Omega\rightarrow \C^N$. 
Notice that  we have an action of $\xi$ on $\O(\Omega)$: given any function $f\in \O(\Omega)$, by normality the function $\xi. f=\L_{\xi}f$ on $\Omega\cap C(Y^{reg})$ extends to a function on $\Omega$. In particular $\xi$ acts on the coordinate functions, so we obtain holomorphic  functions $f_i=\xi. z_i$ on $\Omega$. By possibly making $\Omega$ smaller, we may assume there is a neighborhood $U$ of $\Omega$ in $\C^N$ such that each $f_i$ extends to a holomorphic function on $U$. In particular the vector field $\sum_i f_i \p_{z_i}$ is a holomorphic vector field on $U$ which restricts to $\xi$ on $\Omega\cap C(Y^{reg})$.  For simplicity of notation we also denote  $\xi=\sum_if_i\p_{z_i}$. Now we choose a smaller open set $V\subset\subset U$, then $\xi$ generates a family of local holomorphic transformations $\phi_t$ ($t\in[-\epsilon, \epsilon]$) so that $\phi_t(V)\subset U$. We claim $\phi_t$ maps $\Omega\cap V$ into $\Omega\cap U$. Indeed, given any holomorphic function $f$ on $U$ that vanishes on $\Omega\cap U$,  since $\xi$ is tangent to $\Omega\cap C(Y^{reg})$, we have $\xi. f=0$ on $ \Omega\cap U\cap C(Y^{reg})$ and thus $\xi. f=0$ on $\Omega\cap U$.  This implies that $\phi_t^*f$ vanishes on $\Omega\cap U$ for all $t$, i.e. $f$ vanishes on $\phi_t(\Omega\cap V)$, so the claim follows. 
Clearly $\phi_t$ fixes the vertex $O$, and preserves the function $r$. Using the cone structure it is easy to see that  these local transformations glue together to form  a family of global holomorphic transformations  $\{\phi_t\}_{t\in \bR}$ of $C(Y)$.

 It is also obvious that $\phi_t$ preserves both the smooth and singular part of $C(Y)$.  In particular, it preserves the length of any smooth curve in $C(Y^{reg})$. Using the fact that $C(Y)$ is the metric completion of the Riemannian manifold $C(Y^{reg})$ (Theorem 3.7 in \cite{CC}), it follows that $\phi_t$ acts by isometries on $C(Y)$, hence also on $Y$.  Since $Y$ is compact, by taking the closure of the one-parameter subgroup $\phi_t$ in the isometry group of $Y$ (which is known to be a Lie group by \cite{CC}), we obtain an action of a compact torus $\mathbb{T}$ on $C(Y)$.  This proves Lemma \ref{compact torus action}.\\

The algebraicity of $C(Y)$ depends crucially on this $\mathbb T$ action. Suppose $\Omega$ is a $\mathbb{T}$-invariant neighborhood of $O$. Then we have a weight space expansion
 $$\O(\Omega)=\widehat{\bigoplus}_{\alpha\in \Gamma^*} \O_\alpha(\Omega),  $$
 where $\Gamma^*\subset Lie(\mathbb{T})^*$ is the weight lattice of $\mathbb{T}$, and for $f\in \O_\alpha(\Omega)$, we have $e^{it}. f=e^{i\langle \alpha, t \rangle} f$. The notation $\widehat\bigoplus$ should be understood in terms of Fourier series expansion.  We can take the usual Fourier series expansion of $f$ restricted to each orbit of $\mathbb{T}$, namely, 
 given $f\in \O(\Omega)$, we define 
 \begin{equation}\label{weight decomposition}
 f_{\alpha}(x)=\int_{\mathbb{T}} e^{-i\langle t, \alpha\rangle}f(e^{it}. x)dt. 
 \end{equation}
 It is clear that $f_\alpha$ is holomorphic on $\Omega\cap C(Y^{reg})$, so by normality and continuity $f_\alpha\in \O(\Omega)$.  
Notice the $\mathbb{T}$ action is smooth on $C(Y^{reg})$, so it is easy to see $\sum_\alpha f_\alpha$ converges uniformly on compact subsets of $\Omega\cap C(Y^{reg})$. On the other hand, a singular point  of $C(Y)$ lies in a holomorphic disc with boundary a fixed distance away from the singular set, so by a simple maximal modulus theorem we see the convergence is also uniform on compact subsets of $\Omega$.  

\

Let $N$ be the \emph{embedding dimension} of $C(Y)$ at $O$. This is by definition the smallest integer such that a neighborhood of $O$ embeds holomorphically into $\C^N$.  

\begin{lem} \label{lem3-2}
There is a holomorphic embedding  $F: C(Y)\rightarrow \C^N$ such that the action of $\mathbb T$ extends to a diagonal action on $\C^N$.
\end{lem}

 Choose a local holomorphic embedding $F: \Omega\rightarrow \C^N$ such that $F(O)=0$. By general theory (see for example \cite{GR}, P114-115), $N=\dim_{\C}\m_O/\m_O^2$, where $\m_O$ is the maximal ideal in $\O_O$, 
  and any holomorphic function vanishing on $F(\Omega)$ must have vanishing differential at $0$.

We apply the above expansion to the coordinate functions $z_i=\sum_{\alpha} z_{i, \alpha}$. By Proposition V.B.3 in \cite{GR} there is a polydisc $\Delta$ around $0$, such that each $z_{i, \alpha}$ extends to a holomorphic function on $\Delta$ with a bound $||z_{i,\alpha}||_{L^\infty(\Delta)}\leq C ||z_{i, \alpha}||_{L^\infty(\Omega)}$ for a constant $C>0$ independent of $\alpha$. Therefore we  may assume the series $\sum_\alpha z_{i, \alpha}$ also converges to $z_1$ on $\Delta$.  So there is some $\alpha_i$ such that $dz_{i, \alpha_i}(\p_{z_i})$ is non-zero at $0$. Then the implicit function theorem implies that $F'=(z_{1, \alpha_1}, \cdots, z_{N, \alpha_N})$ is a holomorphic embedding on a possibly smaller neighborhood $\Omega'$ of $O$. Moreover  $F'$ is $\mathbb{T}$-equivariant, where the $\mathbb{T}$ action on $\C^{N}$ is diagonal, with weight on each coordinate given by $\alpha_i$. For simplicity of notation we still denote $(F', \Omega')$ by $(F, \Omega)$.
    
  To extend this  to a global embedding of $C(Y)$, we first notice that since $\xi$ is holomorphic the action of $\mathbb{T}$ induces a holomorphic action of the complexified torus $\mathbb T^\C$ on $C(Y)$. One can see this by first complexifying the action of $Lie(\mathbb{T})$ and then argue as before. Since $F$ is holomorphic and $\mathbb{T}$-equivariant, it is also $\mathbb T^\C$ equivariant, in the sense that if $z$ and $\lambda.z$ are both in $\Omega$ then $F(\lambda.z)=\lambda. F(z)$.  Now we  simply define $F(\lambda. z)=\lambda. F(z)$ for $z\in \Omega$ and $\lambda\in T^{\C}$. Since the radial vector field $r{\p r}=-J\xi$ lies in the Lie algebra of $\mathbb T^\C$, we see that $F$ is defined on $C(Y)$, and it is clear that $F$ is holomorphic. Since $F$ is an embedding near $O$ and $F$ is $\mathbb T^\C$ equivariant, it is also an embedding on the whole  $C(Y)$.  This finishes the proof of Lemma \ref{lem3-2}.
  
 \ 
 
 For simplicity we will call the map $F$ satisfying the property of Lemma \ref{lem3-2} an \emph{equivariant holomorphic embedding}. 
 Now let $W$ be the image of $C(Y)$, endowed with the structure of a reduced complex analytic space. 

\begin{lem}
$W$ is an affine variety in $\C^N$. 
\end{lem}
 Denote by $\mathcal I_W$ the ideal sheaf of $W$. For any $f\in \mathcal I_{W, 0}$, we have a similar expansion $f=\sum_\alpha f_\alpha$ with respect to the $\mathbb{T}$ action on $\C^N$. By the equivariancy $f_\alpha$ also vanishes on $W\cap B$, so  $f_\alpha\in\mathcal I_{W, 0}$. Each $f_\alpha$ extends by homogeneity to an entire holomorphic function on $\C^{N}$ with polynomial growth at infinity, so it must be a homogeneous polynomial.  Therefore  $\mathcal I_{W, 0}$ is generated by the germs of certain homogeneous polynomials. Since $\mathcal I_{W, 0}$ is Noetherian,  it is indeed generated by finitely many of them, say $f_1, \cdots, f_r$.  So $W$ agrees with the affine subvariety in $\C^{N}$ defined by $f_1, \cdots f_r$ in a neighborhood of $0$. By homogeneity they agree globally. 

\

Under the above embedding, the Reeb vector field has an extension to $\C^N$ of the form $\xi=Re(i\sum_{a=1}^Nw_a z_a\p_{z_a})$ for some real numbers $w_1, \cdots, w_N$. 

\begin{lem}
For all $a$, $w_a>0$. 
\end{lem}

For any non-zero polynomial function $f$ on $W$ of weight $\alpha$,  we have $f(\lambda. x)=\lambda^{\langle\alpha, \xi\rangle}f(x)$, where $\lambda.x$ is the radial dilation by $\lambda$ of $x$.  Since $f$ is holomorphic at $0$, it follows that $\langle\alpha, \xi\rangle\geq 0$, and the equality holds if and only if $f$ is radially invariant, i.e. $f$ is a constant. Therefore we see if $\alpha\neq 0$ then $\langle \alpha, \xi\rangle>0$. The lemma follows by applying this to the coordinate functions. \\

Now we describe the affine algebraic structure on $C(Y)$ intrinsically. 
Let $R$ be the ring of holomorphic functions on $C(Y)$ with at most polynomial growth at infinity, and $\H_\alpha$ be the space of polynomial functions on $W$ with weight $\alpha$.  By the above discussion any function $f\in R$ has a Fourier expansion $f=\sum_{\alpha \in \Gamma^*} f_\alpha$, where $f_\alpha\in \H_\alpha$,  and the series converges locally uniformly.  From the formula (\ref{weight decomposition}) each $f_\alpha$ is also of polynomial growth with order at most the growth order of $f$. This implies $\langle\alpha, \xi\rangle$ is uniformly bounded for all $\alpha$ with $f_\alpha\neq 0$. It then follows from the above lemma that there are only finitely many non-zero terms appearing in the expansion. Hence we have a direct sum decomposition
$$R=\bigoplus_{\alpha\in \Gamma^*} \H_\alpha. $$ 
 It is  then straightforward to check that under the above embedding of $C(Y)$ as an affine  variety in $\C^N$, $R$ is naturally identified with the coordinate ring of $W$. In particular, $R$ is finitely generated and $W$ is isomorphic to $\Spec R$. This finishes the proof of Theorem \ref{thm1-2}.

\

 There is also an algebraic description of the cone structure on $C(Y)$. A holomorphic function $f$ on $C(Y)$ is called \emph{homogeneous} with charge $\mu$ if $\L_\xi f=i\mu f$ for some $\mu>0$ (the name is adopted from \cite{GMSY}).  Let $R_d$ be the space of holomorphic functions on $C(Y)$ with charge $d$.  So we may understand the cone structure as a ``grading" on $R$ in terms of the charge: 
 $$R=\bigoplus_{d\in \mathcal S} R_d, $$
where $\mathcal S \subset \bR_{\geq 0}$ is  the \emph{holomorphic spectrum} of $\mathcal S$. Notice the linear map on $\Gamma^*$ sending $\alpha$ to $\langle \alpha, \xi\rangle$ is injective, so each non-zero $R_d$ corresponds to a unique $\H_\alpha$ with $\langle\alpha, \xi\rangle=d$, and we can recover the $\mathbb T$ action from this grading.  We also call the function $h: \mathcal S\rightarrow \Z; d\mapsto \dim R_d$ the \emph{Hilbert function} of $C(Y)$.

 The grading is positive, in the sense  $\xi$ lies in the \emph{Reeb cone} \cite{HS}, \cite{CS}, i.e. the convex cone in $Lie(\mathbb{T})$ consisting of elements $\gamma$ with $\langle\alpha, \gamma\rangle>0$ for all $\alpha\in \Gamma^*$ and $\H_\alpha\neq0$.   Following the terminology introduced in \cite{CS},  we call such $(C(Y), \xi)$ a \emph{polarized affine variety}. 

\

The next result is crucial for the discussion in Section 3.   Notice the Lie algebra $Lie(\mathbb T)$ has a natural rational structure determined by the weight lattice $\Gamma^*$. 

\begin{prop} \label{prop3-5}
The Reeb vector field $\xi\in Lie(\mathbb T)$ is an algebraic vector, as an isolated zero of a system of polynomial equations with rational coefficients. In particular, $\mathcal S$ is contained in the set of algebraic numbers. 
\end{prop}

This is an extension of a result of Martelli-Sparks-Yau \cite{MSY} on the volume minimization property of smooth Sasaki-Einstein metrics. In our setting, the tangent cone  $C(Y)$ admits a (weak) Ricci-flat K\"ahler cone metric, with a global potential given by $r^2$. This enables us to adapt the pluripotential theoretic techniques,  and the proof will be given in the appendix. 

\section{Algebro-geometric description of tangent cones}
\subsection{Rigidity of the holomorphic spectrum}

Let $(Z, p)$ be the Gromov-Hausdorff limit of a sequence of spaces in $\K(n, \kappa, V)$. Recall we have defined a tangent cone at $p$ to be a Gromov-Hausdorff limit of a convergent subsequence of the re-scalings of $(Z, p)$ by $\sqrt{a}$, for integers $a\rightarrow\infty$. 

\begin{lem} \label{lem3-1}
Let $C(Y)$ and $C(Y')$ be two tangent cones at $p$ defined by two sequences of positive integers $\{a_k\}, \{b_k\}$ respectively. Suppose there is a constant $C>0$ so that $C^{-1} \leq a_k/b_k\leq C$ for all $k$,  then $C(Y)$ and $C(Y')$ are isomorphic as elements of $\overline{\K(n, \kappa)}$. In particular, they are isomorphic  as affine algebraic varieties endowed with a Ricci-flat K\"ahler cone metric. 
\end{lem}

This follows from the property of metric cones and the fact from Section 2.3 that the radial dilation on a tangent cone is a holomorphic transformation. 

\

Fix  $\lambda=1/\sqrt{2}$.  Let $(Z_i, p_i)$ be the rescaling of $(Z, p)$ by a factor $\lambda^{-i}$, and we denote by  $\mathcal C_p$ the set of all sequential Gromov-Hausdorff limits of $(Z_i, p_i)$ as $i\rightarrow\infty$. It follows from Lemma \ref{lem3-1} that any tangent cone is indeed isomorphic to one in $\mathcal C_p$.   Notice in Riemannian geometry, the \emph{metric tangent cones} are defined in terms of rescalings of $(Z, p)$ by real numbers $\zeta\rightarrow\infty$ which are not necessarily of the above form $\sqrt{a}$. But a similar argument as Lemma \ref{lem3-1} shows that any general metric tangent cone is also isometric to one in $\mathcal C_p$. So in our context we shall simply call $\mathcal C_p$ the set of tangent cones at $p$. It is endowed with the Gromov-Hausdorff topology. 

\begin{lem} \label{lem3-3}
 $\mathcal C_p$ is compact and connected.
\end{lem}
 
 This should be well-known to experts, and we include a short proof here for the convenience of readers. The compactness follows from Lemma \ref{lem2-0} and the easy fact that $\mathcal C_p$ is a closed subset of $\overline{\K(n, \kappa)}$. Now suppose  $\mathcal C_p$ is a disjoint union of two closed subsets $A$ and $B$. Since $\overline{\K(n, \kappa)}$ is compact and Hausdorff, we can find disjoint open subsets $\mathcal U$, $\mathcal V$ in $\overline{\K(n, \kappa)}$ such that $A\subset \mathcal U$ and $B\subset \mathcal V$. Then it follows that for $i$ sufficiently large $(Z_i, p_i)\in \mathcal U\cup\mathcal V$. Without loss of generality we may assume there is a subsequence $\{\alpha\}\subset \{i\}$ such that $(Z_\alpha, p_\alpha)\in \mathcal U$. Now we claim $(Z_i, p_i)\in \mathcal U$ for all big $i$. For otherwise we may find a subsequence $\{\beta\}\subset\{i\}$ such that $(Z_\beta, p_\beta)\in \mathcal U$, but $(Z_{\beta+1}, p_{\beta+1})\in \mathcal V$. Passing to a subsequence we can assume $(Z_\beta, p_\beta)$ converges to some limit $C(Y)\in \mathcal C_p\cap \overline {\mathcal U}$. Since $\mathcal U$ and $\mathcal V$ are disjoint, and $\mathcal V$ is open, it follows that $\overline {\mathcal U}$ does not intersect $\mathcal V$. Thus we know $C(Y)\in A$. On the other hand,  by Lemma \ref{lem3-1}, $(Z_{\beta+1}, p_{\beta+1})$ also converges to the same limit $C(Y)$, so in particular, for $\beta$ sufficiently large $(Z_{\beta+1}, p_{\beta+1})\in \mathcal U$. Contradiction. Now it follows from the claim that $\mathcal C_p=A$. Hence $\mathcal C_p$ is connected. 
 
 \
 
 Given  a tangent cone $C(Y)\in \mathcal C_p$, the $L^2$ metric over the ball $\{r\leq 1\}$ defines a Hermitian inner product on $R(C(Y))$, which is invariant under the action of $\mathbb T$. Moreover, for $d_1\neq d_2$, $R_{d_1}(C(Y))$ and $R_{d_2}(C(Y))$ are orthogonal.  We have

\begin{thm} \label{thm4-1}
The holomorphic spectrum $\mathcal S:=\mathcal S(C(Y))$ and the Hilbert function of $C(Y)$ are independent of the  tangent cones in $\mathcal C_p$. 
 \end{thm}

For any $D\in \bR^{+}\setminus \mathcal S(C(Y))$, we denote $E_{D}(C(Y))=\bigoplus_{0<d< D}R_d(C(Y))$. 

\begin{lem} \label{lem4-2}
Given $C(Y)$ in $\mathcal C_p$, for any $D\notin \mathcal S(C(Y))$,  there is a small neighborhood $\mathcal U$ of $C(Y)$ so that the vector spaces  $E_{D}(C(Y'))$ have the same dimension for all $C(Y')$ in $\mathcal U$.  
\end{lem}

 By the compactness of $\mathcal C_p$, it suffices to show that if a sequence  $C(Y_j)$  converges to $C(Y)$, then for $j$ sufficiently large $\dim E_j=\dim E$, where we denote $E_j=E_{D}(C(Y_j))$ and $E=E_{D}(C(Y))$.
 
First assume we are given a sequence of homogeneous holomorphic functions $f_j$ on $C(Y_j)$ with charge $d_j\in (0, D)$, and with $||f_j||_{L^2(B_j)}=1$, where $B_j$ the unit ball around the vertex in $C(Y_j)$. Using the interior gradient estimate in Lemma \ref{prop2-10} we obtain a uniform bound of $|\nabla f_j|$ over the half ball $\frac{1}{2}B_j$. By homogeneity for any fixed $k$  we then obtain a uniform bound of $|\nabla f_j|$ over the ball $kB_j$. So by passing to a subsequence $f_j$ converges locally uniformly to a limit $f$ on $C(Y)$. It is clear that $f$ is homogeneous of charge $d\in (0, D)$, and by Lemma \ref{lem2-15} we have $||f||_{L^2(B)}=1$. Now we can apply this to an orthonormal basis of $E_j$. So we conclude that for $j$ big, $\dim E_j\leq \dim E$.  

To prove the other inequality we proceed by contradiction. Suppose $\dim E>\dim E_j$ for all large $j$.  From the above argument, by passing to a subsequence we may assume an orthonormal basis of $E_j$ converges to an orthonormal basis of a proper subspace $E'$ of $E$. Now let $f$ be a function in $E$ which is $L^2$ orthogonal to $E'$, and with $||f||_{L^2(B)}=1$. Suppose $f$ has charge $D_0\in (0, D)$. By Proposition \ref{prop2-8} and using the homogeneity of $f$ we may find for $j$ large a holomorphic function $f_j$ defined on the unit ball $B_j\subset C(Y_j)$ that converges to $f$ uniformly over $B$.  Now using the weight expansion we may write $f_j=g_j+h_j$, where $g_j\in E_j$ and $h_j$ is $L^2$ orthogonal to $g_j$. Then  $||g_j||_{L^2(B_j)}$ is uniformly bounded, so by homogeneity $||g_j||_{L^2(2B_j)}$ is also uniformly bounded.  Using Lemma \ref{prop2-10} again, by passing to a subsequence we may assume $g_j$  converges uniformly to a limit $g$. Hence $h_j$ converges uniformly to $h$, and  $f=g+h$.  By our choice of $f$ we see $g=0$. Now  using the weight expansion for $h_j$ it is easy to see that there is a constant $d\geq D$ such that for all $j$,  $||h_j||_{L^2(\frac{1}{2}B_j)}\leq 2^{-d-n/2} ||h_j|||_{L^2(B_j)}$. Taking limits, this implies $D_0\geq D$. Contradiction.

\

\begin{lem}
There is a dense subset $\mathcal I$ of $\bR^{+}$ such that if $D\in \mathcal I$, then the dimension of $N_D:=E_D(C(Y))$ is independent of $C(Y)\in \mathcal C_p$.
\end{lem}

This follows from the compactness of $\mathcal C_p$ and the fact that $\mathcal S(C(Y))$ is discrete. 

\

Now choose $D\in \mathcal I$. For any $C(Y)\in \mathcal C_p$ we may arrange the holomorphic spectrum of $C(Y)$ within the interval $(0, D)$ (with multiplicities) in the increasing order  as $w_1\leq\cdots \leq w_{N_D}$. From the proof of Lemma \ref{lem4-2} it follows that the map $\iota_D: \mathcal C_p\rightarrow (\bR^+)^{N_D}$ sending $C(Y)$ to its charge vector $(w_1, \cdots, w_{N_D})$ is continuous. 
Since $\mathcal C_p$ is connected, so is the image of $\iota_D$. On the other hand, Proposition \ref{prop3-5} implies that the image is contained in a countable subset of $\bR^{N_D}$, hence it must consist of a single point.  Applying this to all $D\in \mathcal I$, we conclude that $\mathcal S:=\mathcal S(C(Y))$ is independent of $C(Y)$. Then by Lemma \ref{lem4-2} for each $d\in \mathcal S$, $\dim R_d(C(Y))$ is also independent of $C(Y)$. This finishes the proof of Theorem \ref{thm4-1}. 

\

\subsection{Vanishing order of holomorphic functions}

We first set up some notations. Given a tangent cone $C(Y)\in \mathcal C_p$, we denote by $\Lambda$ the dilation by $\lambda$ on $C(Y)$.  Given a function $f$ defined on a ball $B$ in $C(Y)$ around the vertex, we denote by $\Lambda. f$ the function on $B$ with $\Lambda. f(x)=f(\Lambda. x)$.  Let $B_i$ be the unit ball in $Z_i$ around $p_i$. By definition we may naturally identify $B_i$ with a ball in $Z$, and we have natural inclusion maps $\Lambda_i: B_i\rightarrow B_{i-1}$.  For the clarification of later arguments, given a function $f$ defined on $B_{i-1}$, we also denote by $\Lambda_i. f$ the induced function on $B_i$. 
 As $i\rightarrow\infty$, $\Lambda_i$ \emph{converges by sequence} to the dilation $\Lambda: r\mapsto \lambda r$ on the tangent cones\footnote{In this paper, when we say ``converges by sequence", we mean that given any subsequence there is always a further subsequence that converges to some limit. }. There is an ambiguity caused by the possible holomorphic isometric transformation of the tangent cones that fixes the vertex, but this will not affect our following discussion (see Remark \ref{rmk2-10}).

Given a function $f$ defined over a domain in $Z$ that contains $B_i$, we denote by $||f||_{i}$ the $L^2$ norm of the induced function on  $B_i$. If $||f||_i$ is finite, then we define a function $[f]_i$ on $B_i$, which is equal to $a\cdot f|_{B_i}$ for some $a>0$ so that  $||[f]_i||_i=1$.

\begin{lem} \label{lem4-4}
Let $B$ be the ball $\{r<1\}$ in some tangent cone $C(Y)\in \mathcal C_p$.  For any holomorphic function $f$ in $L^2(B)$
 we have 
 $$||\Lambda. f||_{L^2(B)}^2\leq ||f||_{L^2(B)}||\Lambda^2. f||_{L^2(B)}, $$
 and the equality holds if and only if $f$ is homogeneous. 
\end{lem}

We write $f=\sum_{d \in \mathcal S} f_d$, where $f_d$ has charge $d$.  Then $\Lambda. f=\sum_{d\in \mathcal S} \lambda^d f_d$, and $\Lambda^2. f=\sum_{d\in \mathcal S}\lambda^{2d} f_d$. Notice that if $d_1\neq d_2$, then $f_{d_1}$ and $f_{d_2}$ are orthogonal in $L^2(B)$. It follows from the Cauchy-Schwarz inequality that $$||\Lambda. f||_{L^2(B)}^2\leq ||f||_{L^2(B)}||\Lambda^2. f||_{L^2(B)}, $$ and the equality holds if and only if  $f=f_d$ for some $d\in \mathcal S$.  \\

\begin{prop}\label{prop4-5}
For any given $\bar d\notin \mathcal S$, we can find $i_0=i_0(\bar d)$ such that for all $j>i\geq i_0$ and any non-zero holomorphic function $f$ defined on $B_i$,  if $||f||_{i+1}\geq \lambda^{\bar d} ||f||_{i}$, then 
$||f||_{j+1}>\lambda^{\bar d} ||f||_{j}$. 
\end{prop}

 Suppose the conclusion fails, then we would find a subsequence $\{\alpha\}\subset \{i\}$, and non-zero  holomorphic functions $f_\alpha$ defined on $B_{\alpha}$
with 
$$||f_\alpha||_{\alpha+1}\geq \lambda^{\bar d} ||f_\alpha||_{\alpha}$$
$$||f_\alpha||_{\alpha+2}\leq \lambda^{\bar d} ||f_\alpha||_{\alpha+1}$$
 By passing to a subsequence we may assume $B_{\alpha}$ converges to a unit ball $B_\infty$ in some tangent cone. Multiplying $f_\alpha$ by a constant we may assume $||f_\alpha||_{\alpha+1}=1$. Then $||f_\alpha||_{\alpha}\leq \lambda^{-\bar d}$. The gradient estimate Lemma \ref{prop2-10} ensures that by passing to a further subsequence we may assume $f_{\alpha}$ converges to a limit $F$ on $B_\infty$, uniformly on $B_\infty(r)$ for $r<1$. In particular we have 
 $$||F||_{L^2(B_\infty)}\leq \lambda^{-\bar d}, ||\Lambda. F||_{L^2(B_\infty)}=1, ||\Lambda^2. F||_{L^2(B_\infty)}=\lim_{\alpha\rightarrow\infty}||f_{\alpha}||_{\alpha+2}\geq \lambda^{\bar d}. $$
  By Lemma \ref{lem4-4}  $F$ must be  homogeneous holomorphic function on $B_\infty$, and it is clear that the charge must be exactly $\bar d$. This contradicts our hypothesis on $\bar d$.

\begin{cor} \label{cor4-6}
Given a non-zero holomorphic function $f$ defined in a neighborhood of $p\in Z$, then
\begin{enumerate}[(1)]
\item The limit 
$$\lim_{i\rightarrow\infty} (\log\lambda)^{-1}\log( ||f||_{i+1}/||f||_{i})$$
 is either $+\infty$, or a well-defined number in $\mathcal S$. We denote this by $d(f)\in \mathcal S\cup \{+\infty\}$;
\item If $d(f)=+\infty$, then $[f]_i$ converges weakly by sequence to zero;
\item If $d(f)\in \mathcal S$, then $[f]_i$ converges strongly  by sequence to non-zero homogeneous holomorphic functions of charge $d(f)$, on the tangent cones. 

\end{enumerate}

\end{cor}

The existence of $d(f)$ follows immediately from the previous proposition. If $d(f)=\infty$, then by definition all the weak limits must be zero. If  $d(f)\in \mathcal S$, then 
 Lemma \ref{prop2-10} implies that for $i$ large,  $|\nabla f|_{L^\infty(B_i)}\leq K ||f||_{i-1}\leq K\lambda^{-2d(f)} ||f||_i$ for some constant $K>0$ depending only on $n$ and $\kappa$. It then follows that $[f]_i$ converges strongly by sequence. Similar to the proof of the above proposition, any such limit $F$ must satisfy $\lambda^{-2d(f)}||\Lambda^2. F||_{L^2(B_\infty)}=\lambda^{-d(f)}||\Lambda. F||_{L^2(B_\infty)}=||F||_{L^2(B_\infty)}$. 
Again by Lemma \ref{lem4-4}, $F$ must be homogeneous of charge $d(f)$. \\

\begin{rmk}
The above arguments should be compared with classical monotonicity formulas for elliptic equations over cones. The difference is that in our situation we are not exactly working on cones, and this is the place we need to use the rigidity of  the holomorphic spectrum $\mathcal S$.
\end{rmk}

 Notice at this stage we can not rule out the case $d(f)=+\infty$. But later we will do this after establishing the relation with algebraic geometry, see Remark \ref{rmk3-21}. 
Using the estimate in Lemma \ref{prop2-10}, it is easy to see that we have other characterizations
 \begin{eqnarray}\label{eqn4-1}
d(f)&=&\lim_{r\rightarrow0} (\log r)^{-1}\log \sup_{B_r(p)}|f(x)|   \nonumber \\
&=&   \lim_{r\rightarrow0} (\log r)^{-1}\log \sup_{\p B_r(p)}|f(x)|     
\end{eqnarray}

\noindent Hence the number $d(f)$ can be viewed as the vanishing order of $f$ at $p\in Z$, measured by the limit K\"ahler-Einstein metric.  We will use it to study the algebraic geometry of tangent cones. For this purpose we need an extension of Proposition \ref{prop4-5} and Corollary \ref{cor4-6}.  \\

\textbf{Definition.} Suppose we are given a finite dimensional space $P$ of holomorphic functions defined on a neighborhood of $p$. Let $m=\dim P$. 
An \emph{adapted sequence of bases} consists of a basis $\{G_i^1, \cdots, G_i^m\}$ of $P$ for all large $i$,  such that the following holds
\begin{itemize}
\item For all $a$, $||G_i^a||_i=1$; if $a \neq b$, then $\lim_{i\rightarrow \infty} \int_{B_i} G_i^{a}\overline{G_i^{b}}=0$;
\item For  $i$ large and for all  $a$,   $\Lambda_i. G^a_{i-1}=\mu_{ia} G^a_i+p_{i}^a$ for $\mu_{ia}\in \C$, and $p_i^a$ in the linear span $\C\langle G_i^1, \cdots, G_i^{a-1}\rangle$, with $||p_i^a||_i\rightarrow 0$;
\item There are numbers $d_1, \cdots, d_m\in \mathcal S$ with $d_1\leq d_2\leq \cdots\leq d_m$, such that $\mu_{ia}\rightarrow \lambda^{d_a}$;  Moreover,  $p_i^a\in\C \langle G_i^b|b\leq a, d_b=d_a\rangle$. 
\end{itemize}
By definition for $f\in \C\langle G^b_i|a_1\leq b\leq a_2\rangle$ we have $d(f)\in [d_{a_1}, d_{a_2}]$. \\

Now suppose we are given such a space $P$ with an adapted sequence of bases.  

\begin{lem} \label{lem4-8}
$\{[G^a_i]\}$ converges strongly by sequence to an $L^2$ orthonormal set of homogeneous functions of charge exactly $\{d_a\}$, on the tangent cones.  
\end{lem}

 We prove this by induction. For $a=1$ this is clear  by Corollary \ref{cor4-6}. Now we assume the conclusion is true for all $b\leq a-1$. Suppose for a subsequence  $\{\beta\}\subset\{i\}$, $[G^a_{\beta}]$  converges weakly to a limit  $G$, on $B_\infty\subset C(Y)$.  By passing to a further subsequence we may assume $[G^a_{\beta-1}]$ also converges weakly to a limit $G'$ on $B_\infty$. Then by the second and the third item in the above definition, $\Lambda. G'=\lambda^{d_a}G$.  By Lemma \ref{prop2-10}, $\Lambda_\beta. [G^a_{\beta-1}]$ converges uniformly to $\Lambda. G'$. It then follows that $||G||_{L^2(B_\infty)}=\lambda^{-d_a} ||\Lambda. G'||_{L^2(B_\infty)}=\lim_{i\rightarrow\infty} \lambda^{-d_a}||\Lambda_\beta. [G^a_{\beta-1}]||_{\beta}\geq 1$. So $[G^a_\beta]$ converges strongly to $G$. Similarly one can show $\lambda^{-2d_a}||\Lambda^2. G||_{L^2(B_\infty)}=\lambda^{-d_a}||\Lambda. G||_{L^2(B_\infty)}=1$, hence $G$ must be homogenous of charge $d_a$ by Lemma \ref{lem4-4}. By Lemma \ref{lem2-15}, $G$ is $L^2$ orthogonal to the limits of $[G^1_\beta], \cdots, [G^{a-1}_\beta]$. \\
 
 If we choose another adapted sequence of bases, say $\{H_i^a\}$, then by definitions for each $i$, $\{G_i^a\}$ and $\{H_i^a\}$ differ by an action of an element in $U(m)$. So a simple consequence of the above lemma is that the set with multiplicity $d(P)=\{d_1, \cdots, d_m\}$ is independent of the choice of the adapted sequence of bases. 
 
\begin{prop} \label{prop4-9}
For any $\bar d\notin \mathcal S$, we can find $i_0=i_0(\bar d, P)$ such that for all $j>i\geq i_0$, and any holomorphic function $f$ defined on $B_{i}$, if $f\notin P$ and $||\Pi_{i+1} f||_{i+1}\geq \lambda^{\bar d} ||\Pi_{i} f||_{i}$, then 
$||\Pi_{j+1} f||_{j+1}>\lambda^{\bar d} ||\Pi_j f||_j. $
Here $\Pi_j(f)$ denotes the $L^2$ orthogonal  projection of $f|_{B_j}$ to the orthogonal complement of $P|_{B_j}$. 
\end{prop}

 Suppose not, then we may find a subsequence $\{\beta\}\subset\{i\}$, and holomorphic functions $f_\beta$ on $B_\beta$ with 
$$||\Pi_{\beta+1}f_\beta||_{\beta+1}\geq \lambda^{\bar d} ||\Pi_{\beta} f_\beta||_{\beta}$$
$$||\Pi_{\beta+2} f_\beta||_{\beta+2}\leq \lambda^{\bar d} ||\Pi_{\beta+1}f_\beta||_{\beta+1}, $$
We  can normalize so that $||\Pi_{\beta+1} f_\beta||_{\beta+1}=1$. By passing to a subsequence we may obtain weak limits on the ball $B_\infty$ in some tangent cone $C(Y)$: 
$$F=\lim_{\beta\rightarrow\infty} \Pi_\beta f_\beta, F'=\lim_{\beta\rightarrow\infty} \Pi_{\beta+1} f_\beta, F''=\lim_{\beta\rightarrow\infty} \Pi_{\beta+2} f_\beta, $$
with $||F||_{L^2(B_\infty)}\leq \lambda^{-\bar d}$, $||F'||_{L^2(B_\infty)}\leq 1$, and $||F''||_{L^2(B_\infty)}\leq \lambda^{\bar d}$.  
Now we write the $L^2$ orthogonal decomposition on $B_{\beta+1}$
 $$\Lambda_{\beta+1}. \Pi_{\beta} f_\beta =\Pi_{\beta+1} f_\beta+h_{\beta+1},$$
  where $h_{\beta+1}$ is in $P$.  By Lemma \ref{lem4-8}, and by passing to a subsequence, we may assume that both $\{G^a_\beta\}$ and $\{G^a_{\beta+1}\}$ converge to the same orthonormal set $\{G^1, \cdots, G^m\}$ on $B_\infty$, and that $h_{\beta+1}$ converges strongly to a limit $h\in \C\langle G^1, \cdots, G^m\rangle$.  By Lemma \ref{lem2-15} we know $F$ and $F'$ are both orthogonal to $\C\langle G^1, \cdots, G^m\rangle$, so is $\Lambda. F'$,  by the homogeneity of $G^1, \cdots, G^m$. Since $\Lambda. F=F'+h$, we must have $h=0$, and  $F'=\Lambda. F$ with $||F'||_{L^2(B_\infty)}=||\Lambda. F||_{L^2(B_\infty)}\geq 1$. Similarly $F''=\Lambda. F'$.  Then we obtain a contradiction, as in the proof of Proposition \ref{prop4-5}.

\begin{prop} \label{prop4-10}
Suppose $P$ is given as above. Given a holomorphic function $f$ defined on a neighborhood of $p$. Suppose $f\notin P$, then
\begin{enumerate}
\item The limit
 $$\lim_{i\rightarrow\infty} (\log\lambda)^{-1}\log( ||\Pi_{i+1} f||_{i+1}/||\Pi_i f||_{i})$$
 is either $+\infty$ or a well-defined number in $\mathcal S$. We denote this by $d_P(f)$; 
\item If $d_P(f)\in \mathcal S$, then $\hat P=P\oplus \C\langle f\rangle$ also admits an adapted sequence of bases, with $d(\hat P)=d(P)\cup \{d_P(f)\}$.
\end{enumerate}
    \end{prop}

As before, the existence of $d_P(f)$  follows directly from Proposition \ref{prop4-9}.
Now suppose $d_P(f)\in \mathcal S$. We write $\Lambda_i. \Pi_{i-1}f=\Pi_if+h_i$. 
Then for fixed $\epsilon>0$, for all $i$ large,  $||\Pi_{i-1}f||_{i-1}\leq C\lambda^{-d_{P}(f)-\epsilon}||\Pi_if||_i$.  So passing to a subsequence $\{\beta\}\subset \{i\}$ we can take weak limits $C_\alpha^{-1}\Pi_\alpha f\rightarrow F$, $C_\alpha^{-1}\Pi_{\alpha-1}f \rightarrow F'$ and $C_\alpha^{-1}h_\alpha\rightarrow h$, where $C_\alpha=||\Pi_\alpha f||_\alpha$. Then we have $F=\Lambda. F'+h$. Similar to the proof of Proposition \ref{prop4-9} we know $h=0$ and $F$ is homogenous, with  $||F||_{L^2(B_\infty)}=||\Lambda. F'||_{L^2(B_\infty)}=1$. In particular, this implies that $F_i:=[\Pi_i f]_i$ converges strongly by sequence to a homogeneous holomorphic function of charge $d_P(f)$.

Now we write 
\begin{equation} \label{eqn4-2}
\Lambda_i. F_{i-1}=\gamma_i F_i+\sum_{a=1}^{m} \tau_{ia} G^a_i
\end{equation}
 By the above discussion we know $\gamma_i\rightarrow\lambda^{d(f)}$, and $\tau_{ia}\rightarrow 0$ for all $a$. 

\

\textbf{Claim.} We can find $e_{ia}\rightarrow 0$, for  all $a$ with $d_a\neq d(f)$, such that if we replace $F_i$ by  $F_i+\sum_{a: d_a\neq d(f)} e_{ia} G^a_i$, then we can assume that in (\ref{eqn4-2}) $\tau_{ia}=0$ if $d_a\neq d(f)$. 

\

Given this we let $a_0$ be the biggest integer so that $d_{a_0}\leq d(f)$, then we define $\hat G^a_i=G^a_i$ for $a\leq a_0$, $\hat G^{a_0+1}_i=[F_{i}]_i$, and $\hat G^a_i=G^{a-1}_i$ for $a\geq a_0+1$. Then it is easy to check $\{\hat G^a_i\}$ is an adapted sequence of bases for $\hat P$.\\

To prove the claim we let $b$ be the biggest number so that $d_b\neq d_P(f)$. Replacing $F_i$ by $F_i+e_{ib} G^b_i$, then we obtain the new sequence of coefficients $\tilde\tau_{ib}$. In order that  $\tilde\tau_{ib}$ vanishes for all $i$  we need 
$$\gamma_ie_{ib}= \mu_{ib}e_{i-1, b}+\tau_{ib}$$
By assumption, $\lim_{i\rightarrow\infty} \mu_{ib}=\lambda^{d_b}\neq \lambda^{d_P(f)}$, hence by the  lemma below we can choose the desired sequence $\{e_{ib}\}$ for large $i$,  with $e_{ib}\rightarrow 0$. The claim follows by induction on $b$.

\begin{lem} \label{lem4-11}
Given three sequences of complex numbers  $\gamma_i$, $\mu_i, \tau_i$, with limits $\gamma, \mu, \tau$ respectively. Suppose $\tau=0$, and $|\gamma|\neq |\mu|$,   then we can find a sequence $e_i\rightarrow 0$ such that for all $i$ sufficient large, the equation $\gamma_i e_i-\mu_i e_{i-1}=\tau_i$ holds.
\end{lem}

The proof is elementary. If $\gamma=0$ or $\mu=0$ then $e_i$ is uniquely determined and in this case it is easy to see $e_i\rightarrow0$. So we may assume $\gamma \mu\neq 0$. Since we are only interested in large $i$, without loss of generality we may assume for all $i\geq 0$,  $|\gamma_i|\geq |\gamma|/2$ and $|\mu_i|\geq |\mu|/2$.

\

\textbf{Case I}: $|\gamma|>|\mu|$. Then again without loss of generality we may assume there is a $\delta>0$ such that for all $i\geq 0$,  $|\gamma_i|\geq (1+\delta)|\mu_i|$. Set $e_{0}=0$, and define $e_i$ inductively for $i>0$. Then we have
$$|e_{i}|\leq \sum_{j=1}^{i}(1+\delta)^{j-i}|\gamma_{j}|^{-1}|\tau_{j}|$$
For any fixed $l$, we have
$$|e_{i}|\leq 2 (1+\delta)^{l-i}|\gamma|^{-1} \max_{j\leq l}|\tau_j|+2|\gamma|^{-1} (1+\delta)\delta^{-1}\sup_{j\geq l+1} |\tau_j|. $$
It follows easily from this that $e_i\rightarrow 0$. 

\

\textbf{Case II}: $|\gamma|<|\mu|$. In this case we simply define for each $i\geq 0$, 
$$e_{i}=-\sum_{k=1}^{\infty}  \prod_{j=1}^{k-1}(\frac{\gamma_{i+j}}{\mu_{i+j}} )\mu_{i+k}^{-1}|\tau_{i+k}|. $$
By similar arguments as in the previous case we know for all $i$,  this series is absolutely convergent, and $e_i\rightarrow 0$. It is also direct to check this sequence $\{e_i\}$ satisfies the desired equation. This finishes the proof of Lemma \ref{lem4-11}. 

\subsection{Local tangent cones}
Given a tangent cone $C(Y)\in \mathcal C_p$, we denote by $\Aut(C(Y))$ the group of holomorphic transformations of $C(Y)$ that commute with the $\mathbb T$ action generated by the Reeb vector field. We choose  $C(Y)$ such that the dimension of $\Aut(C(Y))$ is minimal among all the tangent cones in $\mathcal C_p$, and we fix a subsequence $\{\alpha\}\subset\{i\}$ that realizes the convergence to $C(Y)$.   As usual we denote the coordinate ring of $C(Y)$ by $R(C(Y))=\bigoplus_{d\in \mathcal S} R_d(C(Y))$. For simplicity we write elements of $\mathcal S$ in an increasing order as $0=d_0<d_1<\cdots$, and we denote $\mu_k=\dim R_{d_k}(C(Y))$, which by Theorem \ref{thm4-1} is independent of the choice of $C(Y)$.

Denote by $\O_p$ the local ring of holomorphic functions defined in a neighborhood of $p\in Z$.  For $d\in \mathcal S$, we let $I_k$ be the subspace of $\mathcal O_p$ consisting of functions with $d(f)\geq d_k$.  This defines a filtration

\begin{equation}\label{eqn4-3}
\O_p=I_{0}\supset I_{1}\supset I_{2}\supset\cdots.
\end{equation}
By (\ref{eqn4-1})  this is a filtration of ideals in $\O_p$, and it is multiplicative in the sense that $I_{j}I_{k}\subset I_{l}$ whenever $d_l\leq d_j+d_k$. Let $R_p$ be the associated graded ring

$$R_p=\bigoplus_{k\geq 0} I_{k}/I_{{k+1}}.$$

\begin{prop} \label{prop4-12}
For all $k\geq 0$, we can find a decomposition $I_{k}=I_{{k+1}}\oplus J_k$, such that $\dim J_k=\mu_{k}$, and $J_k$ admits an adapted sequence of bases with $d(J_k)=\{d_k\}$.
\end{prop}

The proof is by induction on $k$. We first define $J_0$ to be the space of constant functions. 
Now we assume the conclusion holds for all $j\leq k-1$.  Let $\mathcal J$ be the set of all  finite dimensional subspaces $J\subset I_k$ which satisfy $J\cap I_{k+1}=0$, and which admit an adapted sequence of bases with $d(J)=\{d_k\}$.   It is clear by definition that $\dim J\leq \mu_k$ for all $J\in \mathcal J$. Let $J_k \in \mathcal J$ be a maximal element. Now we prove that $\dim J_k=\mu_k$. Suppose not, then  by passing to a subsequence $\{\beta\}\subset\{\alpha\}$ we obtain an orthonormal limit set of homogeneous holomorphic functions $G^1, \cdots, G^p$ on $B_\infty \subset C(Y)$ of charge $d_k$, with $p<\mu_k$.  Now we pick a function $f$ in $R_{d_k}(C(Y))$ with $||f||_{L^2(B_\infty)}=1$, that is orthogonal to $\C\langle G^1, \cdots, G^p\rangle$.
Then by Proposition \ref{prop2-8} for $\beta$ large we may find a sequence of holomorphic functions $f_\beta$ defined on $B_\beta$ that converges uniformly to $f$ as $\beta\rightarrow\infty$.

  Denote $P=\bigoplus_{j\leq k}J_j$.
  Then by Proposition \ref{prop4-5} and \ref{prop4-9}, we see that for $\beta$ large,  $d(f_\beta)\leq d_k$ and $d_P(f_\beta)\leq d_k$.   Fix $\beta_0$ large and let $F=f_{\beta_0}$. By Proposition \ref{prop4-10} we obtain an adapted sequence of bases on $\hat P=P\bigoplus \C\langle F\rangle$ with $d(\hat P)=d(P)\cup \{d_P(F)\}$.  This implies $d_P(F)=d_k$, for otherwise, by taking limits, we obtain a contradiction with the induction hypothesis that $\dim J_j=\mu_j$ for all $j\leq k-1$. From the definition of adapted sequence of bases, we also have $d(F)=d_k$, i.e.  $F\in I_k$, and moreover, $(J_k\bigoplus \C\langle F\rangle)\cap I_{k+1}=0$. It follows that $J_k\bigoplus \C\langle F \rangle \in \mathcal J$, which is strictly bigger than $J_k$. Contradiction. This proves that $\dim J_k=\mu_{k}$. 

To finish the induction step it suffices to prove that $I_{{k}}=I_{k+1}\bigoplus J_k$.  Given any $f\in I_k$ and sufficiently large $\beta$, by rescaling and by adding some element in $J_k$, we may assume that $||f||_\beta=1$ and $f$ is orthogonal to $J_k$ in $L^2(B_\beta)$. Then by Proposition  \ref{prop4-10} we obtain a sequence of adapted bases on $J_k \bigoplus\C\langle f\rangle$. Similar as above, using the fact that $\dim J_k=\mu_k$, we know $d_{J_k}(f)>d_k$. This implies $f\in I_{k+1}\bigoplus J_k$, and hence finishes the proof of Proposition \ref{prop4-12}.

\

Now we fix $D$ large so that  $R(C(Y))$ is generated by $E_{D}(C(Y))$.  Denote $N=\dim E_D(C(Y))$. An orthonormal basis of $E_D(C(Y))$ defines an equivariant embedding  $\Phi: C(Y)\rightarrow\C^N$.  Let $G_\xi$ be the group of linear transformations of $\C^N$ that commute with the $\mathbb T$ action, and let $K_\xi=G_\xi\cap U(N)$. 

Let $k_0=\max\{k\geq 0|d_k<D\}$, and denote $P=\bigoplus_{0<k\leq k_0}J_k$.  By Proposition \ref{prop4-12} we may fix an adapted sequence of bases of $P$, which defines for $i$ large a holomorphic map $F_i: B_i\rightarrow \C^N$, such that the subsequence $F_\alpha$ converges uniformly to $\Phi$ (up to the $K_\xi$ action). Similar to the proof of Proposition \ref{prop2-11}, we may assume $F_i$ is generically one-to-one for all $i$. 

Let $S_k$  be the space of homogeneous  polynomials on $\C^N$ with weighted degree $d_k$, and let $V_k$ be the kernel of the obvious map $S_k\rightarrow R_{d_k}(C(Y))$.  Fix a splitting $S_k=V_k\bigoplus Q_k$, then we may identify $Q_k$ with $R_{d_k}(C(Y))$.   Let $T_{k, \alpha}$ be the subspace of $\O_p$ consisting of the pull back of functions in $Q_k$ by $F_\alpha$.  

\begin{lem} \label{lem4-13}
 Given any $k$, for $\alpha$ large  we have  $I_k=T_{k, \alpha}\bigoplus I_{k+1}$. 
\end{lem}

By the multiplicative property of the filtration we have $T_{k, \alpha}\subset I_k$. By Proposition \ref{prop4-5},  it is easy to see that for $\alpha$ large  $\dim T_{k, \alpha}\geq \mu_k$ and $T_{k, \alpha}\cap I_{k+1}=0$. On the other hand,  by Proposition \ref{prop4-12} $\dim I_k/I_{k+1}=\mu_{k}$. So the lemma follows.

\

Lemma \ref{lem4-13} implies that the ring $R_p$ is finitely generated by $\bigoplus_{k\leq k_0} I_k/I_{k+1}$. Let $W$ be the affine  variety $\Spec(R_p)$. $R_p$ has the same grading as $R(C(Y))$, so $W$ admits a natural action of $\mathbb T$, with the same Hilbert function as $C(Y)$. The chosen adapted sequence of bases of $P$ over $B_i$ then defines a sequence of equivariant embeddings of $W$ into $\C^N$, and we call the image $W_i$.

By general theory (see for example \cite{HaSt}), there is a multi-graded Hilbert scheme $\Hilb$, which is a projective scheme parametrizing polarized affine schemes  in $\C^N$ invariant under the $\mathbb T$ action and with fixed Hilbert function determined by $\{\mu_k\}$.  Therefore $W_i$ (for all large $i$) and $C(Y)$ define points $[W_i]$ and $[C(Y)]$ in $\Hilb$. The group $G_\xi$ acts naturally on $\Hilb$, so that all $[W_i]$ are in the same $G_\xi$ orbit.

 \begin{prop} \label{prop4-14}
$[W_\alpha]$  converges to $[C(Y)]$ in \emph{$\Hilb$}, up to $K_\xi$ action.
 \end{prop} 
 
  By passing to subsequence and by varying $\Phi$ by an element in $K_\xi$, we may fix the ambiguity of $K_\xi$ action and assume that $F_\alpha$ converges to $\Phi$. Fix an arbitrary metric $||\cdot||_*$ on $S_k$.  Given an element $f\in V_k$, for $\alpha$ large  we write $F_\alpha^*f=g_\alpha+F_\alpha^*h_\alpha$ for $g_\alpha\in I_{k+1}$ and $h_\alpha\in Q_k$. We claim that  $||h_\alpha||_*\rightarrow 0$. For otherwise by rescaling we may assume $||h_\alpha||_*=1$, and $g_\alpha+F_\alpha^*h_\alpha=C_\alpha F_\alpha^*f$ with $C_\alpha$ uniformly bounded. Then passing to a subsequence we may assume $g_\alpha$ and $h_\alpha$ converge uniformly to $g$ and $h$ respectively. They satisfy $g+h=0$ and $||h||_*=1$. In particular, $g$ is a non-zero homogeneous function of charge $d_k$.  By Proposition \ref{prop4-5} this would imply  for $\alpha$ sufficiently large that $d(g_\alpha)\leq d_k$. This is a contradiction.  Now we define $f_\alpha=f-h_\alpha\in S_k$. It satisfies that $F_\alpha^*f_\alpha\in I_{k+1}$, so $f_\alpha$ vanishes on ${W_\alpha}$, and $f_\alpha$ converges to $f$ in $S_k$. Now we do the same for a basis of $V_k$ for all $k\leq k_1$, where $k_1$ is chosen so that any ideal of $\C[x_1, \cdots, x_N]$ defining an element in $\Hilb$ is generated by the homogeneous pieces of degree at most $k_1$. It then follows that $[C(Y)]$ is the limit of $[W_\alpha]$ in $\Hilb$. \\

  Since the universal family over $\Hilb$ is flat and normality is an open condition in a flat family (see for example \cite{Banica}), it follows that $W$ is normal variety. 
 Recall for all $C(Y')\in \mathcal C_p$,  a choice of orthonormal basis of $E_D(C(Y'))$ determines a holomorphic map $\Phi': C(Y')\rightarrow\C^N$. 

\begin{lem} \label{lem4-15}
 There is a neighborhood $\mathcal U$ of $C(Y)$ in $\mathcal C_p$ such that for all $C(Y')\in \U$, $\Phi'(C(Y'))$ is normal.
 \end{lem}
 
Otherwise we choose  a sequence $C(Y_s)$ converging to $C(Y)$ such that the image $\Phi_s(C(Y_s))$ is not normal. By modifying $\Phi_s$ by elements in $K_\xi$, we may assume $\Phi_s(C(Y_s))$ converges to $F(C(Y))$. Now for each $s$, we can find $D_s$ big so that $R(C(Y_s))$ is generated by elements of charge at most $D_s$. Choose a subsequence $\{\beta\}\subset \{i\}$ so that $B_\beta$ converges to the unit ball in $C(Y_s)$. Then we can argue as above replacing $D$  by $D_s$ and $N$ by $N_s=\dim E_{D_s}(C(Y_s))$, and assume $W_\beta$ converges to $C(Y_s)$ as affine varieties in some $\C^{N_s}$, i.e. the convergence is taken in a different multi-graded Hilbert scheme. Projecting down to $\C^N$, we see that $W_\beta$ converges to $F_s(C(Y_s))$ locally as complex analytic spaces. Now since $\Hilb$ is  compact, by passing to a subsequence we may also assume $[W_\beta]$ converges to a limit $[\Sigma_s]$ in $\Hilb$. It then follows that  the underlying reduced complex analytic space of $\Sigma_s$ is the same as $\Phi_s(C(Y_s))$. So $\Sigma_s$ converges to $C(Y)$ locally as complex analytic spaces in $\C^N$. Now using the compactness of $\Hilb$ again by passing to a subsequence we may assume $[\Sigma_s]$ converges to a limit $[\Sigma]$, whose underlying reduced complex analytic space coincides with $C(Y)$. On the other hand, since $C(Y)$ is normal and $[\Sigma]$ and $[C(Y)]$ have the same Hilbert function, it follows that $[\Sigma]=[C(Y)]$. This implies by openness of normality again that $\Sigma_s$ is normal for $s$ large. In particular we know $\Phi_s(C(Y_s))$ is normal. Contradiction.

 \
 
 Now we prove Theorem \ref{thm1-3}. 
By making $\mathcal U$ even smaller, we may assume by Proposition \ref{prop2-11} that $\Phi'$ is generically one-to-one and so by Lemma \ref{lem4-15} $\Phi'$ is an embedding. In particular, $C(Y')$ also defines an element $[C(Y')]$ in $\Hilb$.
 From the construction \cite{HaSt},  $\Hilb$ is  a sub-scheme of a certain projective space $\P$, and the action of $G_\xi$ extends to $\P(V)$. It follows easily from the definition that the stabilizer of $[C(Y)]\in \Hilb$ is isomorphic to $\Aut(C(Y))$ which by Proposition \ref{propA-6},  is reductive. So we can write $\Aut(C(Y))=K^\C$, for a compact group $K$. 
 
  As in \cite{Do10} (Proof of Proposition 1), we can find an equivariant slice for the action. Namely, there is a projective subspace $\P'=\P(\C v\oplus S)$, where $v$ is a vector in $V$ lying over $[C(Y)]$ and $S$ is a $K^\C$ invariant subspace of $V$ which is transverse to the $G_\xi$ orbit of $[C(Y)]$. Let $O$ be the $G_\xi$ orbit of $[W]$, and $O'=O\cap \P'$. Notice by general theory the closure $\overline {O'}$ is a (possibly reducible) algebraic variety. By Proposition \ref{prop4-14} we know $[C(Y)]\in \overline{O}$. So from the construction of $\P'$ in \cite{Do10} we can find a small neighborhood $U$ of $[C(Y)]$ in $\P$, such that  each component of $O'\cap U$  is contained in a single $K^\C$ orbit. Moreover any point in $\overline O\cap U$ is in the $G_\xi$ orbit of a point in $\overline{O'}\cap U$. In particular, $[C(Y)]\in \overline{O'}$.
  
   Suppose $C(Y')\in \mathcal C_p$ is close to $C(Y)$, then we may assume $[C(Y')]\in \overline O \cap U$. By the above discussion we may find $g\in G_\xi$ such that $g. [C(Y')]\in \overline{O'}\cap U$. We claim $[C(Y)]$ is in the closure of the $K^\C$ orbit  of $g. [C(Y')]$.  Indeed, since $[C(Y)]$ is fixed by $K^\C$, we may reduce to the linear action on $S$,  and this becomes the well-known fact that if $x\in S$ is such that $0\in \overline{K^\C. x}$, then for any $y\in \overline{K^\C. x}$, we have $0\in \overline{K^{\C}. y}$ (the point is that $0$ is a \emph{closed} $K^\C$ orbit, and any $K^\C$ invariant polynomial on $S$ vanishing at $y$ must also vanish at $0$). The claim implies that $[C(Y)]$ and $[C(Y')]$ are in the same $G_\xi$ orbit, for otherwise we would have $\dim \Aut(C(Y))>\dim \Aut(C(Y'))$, which contradicts our choice of $C(Y)$. By the uniqueness of Ricci-flat K\"ahler cone metric (Proposition \ref{propA-7}) on $C(Y)$, it follows that $[C(Y)]$ and $[C(Y')]$ are isomorphic as affine varieties with a  Ricci-flat K\"ahler cone metric, and so are indeed in the same $K_\xi$ orbit.  Theorem \ref{thm1-3} then follows from  the connectedness of $\mathcal C_p$.
 
\begin{rmk}\label{rmk3-17}
The precise meaning of Theorem \ref{thm1-1} is that any two tangent cones are isomorphic as affine algebraic varieties endowed with a Ricci-flat K\"ahler cone metric. Notice by the discussion of Section 2.2 each tangent cone is also given a polarization, which is a priori not unique from the definition. It seems an interesting question to further examine the limiting polarization, in particular the $U(1)$ connection. We leave this for future study.
\end{rmk}

A consequence of the above argument, using the Hilbert-Mumford criterion, is that there is a one parameter subgroup $\lambda(t)$ of $G_\xi$, such that $[C(Y)]=\lim_{t\rightarrow0} \lambda(t). [W]$. In terms of the language of K-stability, we may say there is a \emph{test configuration}  for $W$, with central fiber $C(Y)$, in the sense of \cite{CS}. 

\

Now we study the meaning of $W$ in terms of the local complex analytic geometry of $Z$ at  $p$. 
First we recall the notion of a \emph{weighted tangent cone}. Let $(w_1, \cdots, w_m)\in (\bR^+)^m$ be a weight vector, and assign any monomial $z_1^{a_1}\cdots z_m^{a_m}$ with weight $\sum_i a_iw_i$.  For any holomorphic function $f$ defined in a neighborhood of $0\in \C^m$, we let $w(f)$ be the smallest weight among all monomials in the Taylor expansion of $f$. Suppose $(X, 0)$ is a germ of a complex analytic set in $\C^m$.   Consider the weight filtration
$$\O_0=\mathcal F_{e_0}\supset \mathcal F_{e_1}\supset\cdots$$
where $\mathcal F_{e_k}$ consists of the restriction of holomorphic functions $f$ on a neighborhood of $0$ with $w(f)\geq e_k$. The associated graded ring $\mathcal R(\F)=\bigoplus_{k\geq 0}\mathcal F_{e_k}/\mathcal F_{e_{k+1}}$ is naturally isomorphic to $\C[x_1, \cdots, x_m]/\mathcal I$, 
where $\mathcal I$ is the ideal generated by weighted homogeneous polynomial functions $f$ on $\C^m$ such that $f|_X$ is equal to the restriction of a germ of analytic function $g$ with $w(g)>w(f)$, i.e. $f$ is the initial term of a defining equation of $X$ at $0$ (with respect to the above weight). Therefore $\Spec(\R(\F))$ defines a polarized affine sub-scheme in $\C^m$, with Reeb vector field $\xi=\sum_{a}Re(iw_az_a\p_{z_a})$.  We call it the weighted tangent cone of $(X, 0)$.   Notice if all the weights are equal, then $\Spec(\R(\F))$ is  the \emph{Zariski tangent cone} of $X$ at $0$, which is independent of the choice of analytic embedding. In general however, the weighted tangent cone depends on the choice of the analytic embedding, but it is invariant under the action of $G_\xi$. In particular, in our situation above for the obvious weight vector, the weighted tangent cones of $(F_i(B_i), 0)$ are all isomorphic.

\begin{prop} \label{prop4-17}
$W$ is isomorphic to  the weighted tangent cone of $(F_i(B_i), 0)$ in $\C^N$, with respect to the weight determined by the $\mathbb T$ action.
\end{prop}

Without loss of generality we may assume $i=1$.  We define a natural map from $\C[x_1, \cdots, x_N]$ to $R_p$, that sends a polynomial $f$ with $w(f)=d_k$ to the subspace $I_k/I_{k+1}$. This is well-defined since  $d(F_1^*f)\geq w(f)$ by (\ref{eqn4-1}). It also descends to a map $\tau: \R(\F)\rightarrow R_p$. By Lemma \ref{lem4-13} $\tau$ is surjective. So it suffices to show $\tau$ is also injective. For this we need a lemma. Let $\widetilde O$ be the sheaf of holomorphic functions on $F_1(B_1)$. For simplicity of notation we view $\widetilde O_0$ as a subspace of $\O_p$ via the obvious map.

\begin{lem} \label{lem4-18}
There is a function $d'=d'(d)$ that grows linearly as $d\rightarrow\infty$,  such that if a holomorphic function $f\in \widetilde O_0$ satisfies $f\in I_d$, then $f\in \m_0^{d'}$, where $\m_0$ is the maximal ideal in $\widetilde O_0$. 
\end{lem}

 Given this,  suppose $f\in \F_{d_k} \cap I_{k+1}$, then using Lemma \ref{lem4-13}  for $\alpha$ large  we can write
$$f=f_1+f_2+\cdots+f_l+g_l$$ 
with $f_l\in T_{l, \alpha}\subset \F_{d_{k+l}}$ and $g_l\in I_{d_{k+l}}$. By Lemma \ref{lem4-18} and the fact that all the weights are positive we know if we make $l$ sufficiently large, then $g_l\in \F_{d_{k+1}}$. This shows the map $\tau$ is injective and finishes the proof of Proposition \ref{prop4-17}. 
 
 \
 
 It remains to prove Lemma \ref{lem4-18}.
In $\C^N$ we define  $||x||^2=(\sum_{a}|x_a|^{2/w_a})^{1/2}$. 
From the definition of the adapted sequence of bases it it easy to see that 
for any $\epsilon>0$ small, there is a constant $C_\epsilon>0$, such that for all $x\in B_1$
$$C_\epsilon ||F_1(x)||^{1+\epsilon}\leq d_Z(x, p)\leq C_\epsilon ||F_1(x)||^{1-\epsilon}.$$
So if $f\in I_d$ then we have $|f(x)|\leq C_\epsilon' ||x||^{d-2\epsilon d}$ for some constant $C'_\epsilon>0$. Now we first blow up $F(B_1)$ at $0$ and then let $\hat B$ be a resolution of singularities of the blown-up. Let $\pi: \hat B\rightarrow F(B_1)$ be the natural projection map, then by general theory $\pi^{-1}\m_0=\O(-\sum b_i E_i)$ where $E_i$ are the exceptional divisors over $0$, and $b_i$ are positive integers. Clearly on compact sets of $\C^N$,  $||x||$ is H\"older equivalent to the Euclidean norm, so using the above estimate of $f$ we see that $\pi^*f$ has vanishing order at least $Cd$ along each $E_i$ for some constant $C>0$. Hence the lemma follows.

\begin{rmk} \label{rmk3-21}
By the Krull intersection theorem  $\bigcap_{d\geq0} \m_0^d=0$, so using Lemma \ref{lem4-18} we have $\bigcap_{k\geq0} I_k=0$. In particular, this implies that  $d(f)$ is indeed finite for any non-zero function $f\in \mathcal O_p$. 
\end{rmk}

 \

From Proposition \ref{prop4-17} we obtain a flat family of complex analytic spaces with central fiber $W$ and general fiber $F_1(B_1)$.  Using openness of normality again it follows that $F_1(B_1)$ is normal. Since $F_1$ is generically one-to-one, we conclude that $F_1$ is a holomorphic equivalence.

\

 To sum up, we have achieved the following:
 
 \begin{itemize} 
\item There is a unique tangent cone $C(Y)$ of $Z$ at $p$, as an affine algebraic variety together with a Ricci-flat K\"ahler cone metric. 
\item There is a polarized affine variety $W$, obtained as a weighted tangent cone of $Z$ at $p$ under some  local holomorphic embedding;
\item There is a test configuration for $W$ as a polarized affine algebraic variety, with central fiber $C(Y)$.
\end{itemize}

\

\noindent \textbf{Further discussion}:

\

 In algebraic geometry, it is a classical fact that  the Zariski tangent cone is  an intrinsic object associated to a germ of singularity. The above weighted tangent cone $W$ is usually not the same as the Zariski tangent cone, but we expect that both $W$ and $C(Y)$ are also intrinsic invariants of a local algebraic singularity. Notice as in \cite{WN}, \cite{Sz},  a test configuration for $W$ can also be viewed as a filtration on the co-ordinate ring of $W$. In terms of the notion of K-stability for polarized affine varieties formulated in \cite{CS}, and suppose the results of \cite{Berman, CDS1, CDS2, CDS3} extend to this case, we can say $W$ is \emph{K-semistable} and $C(Y)$ is \emph{K-stable}.  So we see some similarity between the above picture and the well-known Harder-Narasimhan filtration/Jordan-H\"older filtration for holomorphic vector bundles.  This also motivates the following 
 
 \begin{conj}
The filtration (\ref{eqn4-3}) and the polarized affine varieties $W$ and $C(Y)$ are uniquely determined by the germ of the analytic singularity $p$. In particular they are independent of the K\"ahler-Einstein metric on $Z$ defining them.  
 \end{conj}
 
One can also formulate a corresponding algebro-geometric conjecture characterizing $W$ and $C(Y)$ in terms of K-stability.

\

It is an interesting question to understand these for explicit  algebraic singularities. Here we discuss two classes of examples. First  we consider a class of isolated hypersurface singularities. For $n\geq 2$ and $k\geq 1$, we denote by $X^n_k$  the hypersurface in $\C^{n+1}$ $(n\geq 2)$ with defining equation $x_0^{k+1}+x_1^2+\cdots +x_n^2=0$.  The origin $0$ is the unique singular point, and the germ at $0$ is  usually called an  \emph{n dimensional $A_k$ singularity}.  We divide the range of $(n, k)$ into three categories: 
$$\I=\{(n,k)|n=2\}\bigcup \{(n, k)|n=3, k\leq 2\}\bigcup \{(n, k)|k=1\};$$
$$\II=\{(n,k)|n=3, k\geq 4\}\bigcup\{(n,k)|n=4, k\geq 3\};$$
$$\III=\{(3, 3)\}\bigcup \{(4, 2)\}.$$
Under the above embedding in $\C^{n+1}$, $X^n_k$ is naturally a polarized affine variety with respect to the obvious weight  $(2, k+1, \cdots, k+1)$. By \cite{GMSY, LS},  there is a compatible Ricci-flat K\"ahler cone metric on $X^n_k$ if and only if $(n,k)$ belongs to $\I$. When $n=2$ this is the flat orbifold cone. When $k=1$ this is the $n$-dimensional Stenzel's cone. 

Now suppose our limit space $(Z, p)$ is locally analytically isomorphic to  $(X^n_k, 0)$. Let $C(Y)$ be the tangent cone at $p$, and $W$ be the affine variety obtained as above. 
The question is to describe $W$ and $C(Y)$. Our conjectural picture depends on the range of  $(n, k)$

\begin{enumerate}
\item $(n,k)\in \I$. Naturally one expects that $W$ and $C(Y)$ are both the known Ricci-flat K\"ahler metric with the standard Reeb vector field.
\item $(n,k)\in \II$. In this case,  Hein-Naber \cite{HN} constructed a Calabi-Yau metric in a neighborhood of $0$ in $X^n_k$, with the tangent cone at $0$ given by $X_\infty^{n}$. Here $X_\infty^n$ is the hypersurface in $\C^{n+1}$ defined by $x_1^2+\cdots+x_n^2=0$, endowed with the product of the $n-1$ dimensional Stenzel cone and the flat metric on $\C$. We expect that in general $W$ and $C(Y)$ are both isomorphic to $X^n_\infty$. 
\item  $(n,k)\in \III$: These are critical cases,  and we expect $W$ is isomorphic to $X^n_k$, but $C(Y)$ is isomorphic to $X^n_\infty$.
\end{enumerate}

We make some simple observations that support this picture. It is not hard to see that under the natural embedding into $\C^{n+1}$, the Ricci-flat K\"ahler cone metric on $X^n_\infty$ has weight vector given by $w=(1, 2\frac{n-1}{n-2}, \cdots, 2\frac{n-1}{n-2})$. If we consider the standard embedding of $X^n_k$ in $\C^{n+1}$, then one sees that the weighted tangent cone with respect to $w$ is given by $X^n_\infty$ exactly when $(n, k)\in \II$. In the case $(n, k)\in \III$, $X^n_k$ is itself a polarized affine variety with respect to $w$, so it is natural to hope that $X^n_k$ degenerates to $X^n_\infty$ by another $\C^*$ action that is equivariant with respect to $w$(which is obvious to find).  With slightly more work, one can show that when $(n, k)\in \I$,  the tangent cone can never be $X^n_\infty$. In general it still remains an algebro-geometric question to verify the above expectations.  We leave this for future work. Notice by Proposition \ref{prop2-13} if $p$ is a smooth point of $Z$ (in the complex-analytic sense), then both $W$ and $C(Y)$ are isomorphic to $\C^n$ (with the standard weight), but to our knowledge even in this case a purely algebro-geometric proof of this fact is still lacking. 

\

For another class of examples, we suppose $(Z, p)$ is toric, i.e. there is an effective action of an $n$-dimensional torus $T^n$ on $Z$ that fixes $p$ and preserves the limit metric and complex structure. This happens when $(Z, p)$ is a  toric $\Q$-Fano variety (by the uniqueness of K\"ahler-Einstein metrics \cite{BBEGZ}). There are interesting examples appearing on the boundary of the compactification of smooth Fano manifolds (c.f. \cite{Spotti}).   for Let $\Delta_Z$ be the moment polytope of $Z$. In this case one can see that the above discussion can be made in a $T^n$-equivariant manner.  In particular, both $W$ and $C(Y)$ are also toric.  Moreover, there is a $T^n$-invariant neighborhood $U$ of $p$, and a holomorphic embedding of $U$ into some $\C^N$ such that the $T^n$ action extends to a diagonal action on $\C^N$,  and $W$ is realized as a weighted tangent cone of $U$ at $p$.  The fact that $W$ is normal implies that the Reeb vector field of $W$ indeed lies in the Lie algebra of $T^n$. In particular, as an affine toric variety, $W$ is isomorphic to the natural toric tangent cone $T_pZ$ of $Z$ at $p$, with moment polytope given by the Euclidean tangent cone of $\Delta_Z$ at $p$.  Similarly one can show that $C(Y)$ is also isomorphic to $W$ as polarized affine varieties.  It seems possible that one can further adapt the results in \cite{MSY2} to this situation, and determine the Reeb vector field of $C(Y)$ inside the Reeb cone, in terms of the geometry of $T_pZ$.

\subsection{Tangent cones at infinity}

Now we turn to tangent cones at infinity. The results will be mostly parallel to the case of local tangent cones. 
Let $(Z, p)$ be a Gromov-Hausdorff limit of a sequence of spaces in $\K(n, \kappa, V)$, and we assume $Z$ is non-compact, i.e. the rescaling factors $a_i\rightarrow\infty$. Again fix $\lambda=1/\sqrt{2}$, and let $(Z_i, p_i)$ be the rescaling of $(Z, p)$ by  $\lambda^{i}$ this time.  A tangent cone at infinity is a Gromov-Hausdorff limit of a convergent subsequence of $(Z_i, p_i)$. It is clear that a tangent cone itself is also a Gromov-Hausdorff limit of spaces in $\K(n, \kappa, V)$  with the rescaling factors tending to infinity.  Let $\mathcal C_\infty$ be the set of all tangent cones at infinity. These are independent of the choice of base point $p$.

It is straightforward to adapt the results of Section 2.3 and Section 3.1 to show that any tangent cone $C(Y)\in \mathcal C_\infty$ is a polarized affine algebraic variety with coordinate ring $R(C(Y))=\bigoplus_{k\geq 0}R_{d_k}(C(Y))$, and the holomorphic spectrum $\mathcal S=\{d_k\}$ is independent of $C(Y)$. We also have analogous results to Section 3.2, with almost identical proofs. For the convenience of readers we write down the statements here, and only point out the part of proof that is different from the case of local tangent cones.  We adapt the notations at the beginning of Section 3.2, except the natural inclusion map is now given by $\Lambda_i: B_i\rightarrow B_{i+1}$. 

\begin{prop} \label{prop4-20}
For any given $\bar d\notin \mathcal S$, we can find $i_0=i_0(\bar d)$ such that for all $j>i\geq i_0$ and any non-zero holomorphic function $f$ defined on $B_j$,  if $||f||_{j}\leq \lambda^{-\bar d} ||f||_{j-1}$, then 
$||f||_{i}<\lambda^{-\bar d} ||f||_{i-1}$. 
\end{prop}

As in (\ref{eqn4-1}), given a holomorphic function $f$ on $Z$, we can define the order of growth at infinity by
\begin{equation}\label{eqn4-4}
d(f)=\lim_{r\rightarrow\infty} (\log r)^{-1}\sup_{B_r(p)} \log |f(x)|. 
\end{equation}
Similar to the proof of Corollary \ref{cor4-6}, this is well-defined and one can show $d(f)\in \mathcal S\cup\{+\infty\}$. Now let $R(Z)$ be the ring of all holomorphic functions $f$ on $Z$ with polynomial growth (i.e. with $d(f)<+\infty$).

\

 Given a finite dimensional subspace $P\subset R(Z)$ with dimension $m$. We can similarly define the notion of an \emph{adapted sequence of bases}. It consists of a basis $\{G_i^1, \cdots, G_i^m\}$ of $P$ for all large $i$,  such that the following holds
\begin{itemize}
\item For all $a$, $||G_i^a||_i=1$; if $a \neq b$, then $\lim_{i\rightarrow \infty} \int_{B_i} G_i^{a}\overline{G_i^{b}}=0$;
\item For $i$ large and all $a$,   $\Lambda_i. G^a_{i+1}=\mu_{ia} G^a_i+p_{i}^a$ for $\mu_{ia}\in \C$, and $p_i^a\in \C\langle G_i^1, \cdots, G_i^{a-1}\rangle$, with $||p_i^a||_i\rightarrow 0$;
\item There are numbers $d_1, \cdots, d_m\in \mathcal S$ with $d_1\leq d_2\leq \cdots\leq d_m$, such that $\mu_{ia}\rightarrow \lambda^{d_a}$;  Moreover,  $p_i^a\in\C \langle G_i^b|b\leq a, d_b=d_a\rangle$. 
\end{itemize}
Again for $f\in \C\langle G^b_i|a_1\leq b\leq a_2\rangle$ we have $d(f)\in [d_{a_1}, d_{a_2}]$. We also define $d(P)=\{d_1, \cdots, d_m\}$.\\

\begin{prop} \label{prop4-21}
For any $\bar d\notin \mathcal S$, we can find $i_0=i_0(\bar d, P)$ such that for all $j>i\geq i_0$, and any holomorphic function $f$ defined on $B_{j}$, if $f\notin P$ and $||\Pi_{j} f||_{j}\leq \lambda^{-\bar d} ||\Pi_{j-1} f||_{j-1}$, then 
$||\Pi_{i} f||_{i}<\lambda^{-\bar d} ||\Pi_{i-1} f||_{i-1}. $
Here $\Pi_j(f)$ denotes the $L^2$ orthogonal  projection of $f|_{B_j}$ to the orthogonal complement of $P|_{B_j}$. 
\end{prop}

\begin{prop} \label{prop4-22}
Given a holomorphic function $f\in R(Z)$. Suppose $f\notin P$, then the following limit
 $$-\lim_{i\rightarrow\infty} (\log\lambda)^{-1}\log( ||\Pi_{i+1} f||_{i+1}/||\Pi_i f||_{i}) $$
 is a well-defined element in $\mathcal S\cup\{+\infty\}$, which we denote by $d_P(f)$. 
  Moreover, if $d_P(f)\in \mathcal S$, then $\hat P=P\oplus \C\langle f\rangle$ also admits an adapted sequence of bases, with $d(\hat P)=d(P)\cup \{d_P(f)\}$.
\end{prop}
    
    Now we fix a tangent cone $C(Y)\in \mathcal C_\infty$, and a subsequence $\{\alpha\}\subset\{i\}$ such that $B_\alpha$ converges to the unit ball $B$ in $C(Y)$. 
  For $d\in \mathcal S$, we denote by $I_d$ the space of holomorphic functions $f$ on $Z$ with $d(f)\leq d$.  Again we list elements in $\mathcal S$ with increasing order $0=d_0<d_1<\cdots$, and denote $\mu_k=\dim R_{d_k}(C(Y))$. Then we define a filtration of  $R(Z)$ 
\begin{equation}\label{eqn4-5}
0=I_{0}\subset I_{1}\subset I_{2}\subset\cdots, 
\end{equation}
and correspondingly a graded ring 
$$R_\infty(Z)=\bigoplus_{k\geq 0}I_{k+1}/I_k.$$
The difference from (\ref{eqn4-3}) is that the inclusion direction is reversed.

\begin{prop}  \label{prop4-23}
For all $k\geq 0$, we can find a decomposition $I_{k+1}=I_{{k}}\oplus J_k$, such that $\dim J_k=\mu_k$, and $J_k$ admits an adapted sequence of bases with $d(J_k)=\{d_k\}$.
\end{prop}

The proof is similar to Proposition \ref{prop4-12}, except a new technical point due to the fact that a priori $I_k$ may be empty for $k\geq 1$ so we need to construct global holomorphic functions in the meantime. Again we prove by induction on $k$. Let $J_0$  be the space of constant functions. 
Now we assume the conclusion holds for all $j\leq k-1$.  Let $\mathcal J$ be the set of all  finite dimensional subspaces $J\subset I_k$ which satisfy $J\cap I_{k-1}=0$, and which admit an adapted sequence of bases with $d(J)=\{d_k\}$.  By definition of adapted sequence of bases we have for all such $J$ that $\dim J\leq \mu_k$. Let $J_k \in \mathcal J$ be a maximal element.  Now we prove that $\dim J_k=\mu_k$. Suppose not, then  by passing to a subsequence $\{\beta\}\subset\{\alpha\}$ we obtain an orthonormal limit set of homogeneous holomorphic functions $G^1, \cdots, G^p$ on $B_\infty \subset C(Y)$ of charge $d_k$, with $p<\mu_k$.  Now we pick a function $f$ in $R_{d_k}(C(Y))$ with $||f||_{L^2(B_\infty)}=1$, that is orthogonal to $\C\langle G^1, \cdots, G^p\rangle$.
Then by Proposition \ref{prop2-8} for $\beta$ large we may find a sequence of holomorphic functions $f_\beta$ defined on $B_\beta$ that converges uniformly to $f$ as $\beta\rightarrow\infty$.

  Denote $P=\bigoplus_{j\leq k}J_j$. By by Proposition \ref{prop4-20} and \ref{prop4-21}, we can find $\epsilon>0$ small and $\beta_0>0$ such that if $\beta>\beta_1>\beta_0$ then $||f_{\beta_1+1}||_{\beta_1+1}\leq \lambda^{-d_k-\epsilon} ||f_{\beta_1}||_{\beta_1}$, and $||\Pi_{\beta_1+1}f_{\beta_1+1}||_{\beta_1+1}\leq \lambda^{-d_k-\epsilon} ||\Pi_{\beta_1}f_{\beta_1}||_{\beta_1}$. So in particular, by using a diagonal sequence argument we may assume $f_\beta$ converges to a limit $F$ over any fixed size ball, with $d(F)\leq d_k$ and $d_P(F)\leq d_k$. In particular $F\in I_{k}$. From this point, the proof proceeds identically the same as Proposition \ref{prop4-12}, and we skip it here.

\

Now we proceed to prove Theorem \ref{thm1-4}. Choose $k_0$ so that $d_{k_0}>D$ and $R(C(Y))$ is generated by $E_D(C(Y))$, and  a sequence of adapted bases for $J_k$ for all $k\leq k_0$. Using these we define  maps $F_i: Z\rightarrow \C^N$ for $i$ sufficiently large with $F_i(p_i)=0$, where $N=\dim E_D(C(Y))$. 

As in Lemma \ref{lem4-13}, one then proves that $R(Z)$ is generated by $\bigoplus_{k\leq k_0}J_k$ and $R_\infty(Z)$ is generated by $\bigoplus_{k\leq k_0}I_{k+1}/I_{k}$.  The chosen adapted bases of $J_k$ for $k\geq k_0$ then realizes $\Spec(R(Z))$ as an affine  variety $\tilde Z_i$ in $\C^N$.  It is clear that $F_i(Z)\subset \tilde Z_i$. We claim that $F_i(Z)=\tilde Z_i$. Notice by definition $R(Z)$ is an integral domain, so $\tilde Z_i$ is reduced and irreducible. Thus it suffices to prove that $\dim \tilde Z_i=\dim Z$.  For this we notice that $\dim I_k=\sum_{j=0}^k R_{d_j}\leq C d_k^{n}$, so the dimension of polynomial functions on $\tilde Z_i$ with the usual degree at most $d$ is also bounded by $Cd^{n}$,  and hence $\dim \tilde Z_i\leq n$.  

\

Now we can follow the same arguments as in Section 3.3 to show further than $F_i$ is indeed a holomorphic embedding, and furthermore, there is a unique tangent cone $C(Y)$ at infinity. This finishes the proof of Theorem \ref{thm1-4}. Moreover,  one can obtain an algebro-geometric description of the tangent cone at infinity, similar to Section 3.3. However, in general one would not expect a naive intrinsic algebro-geometric characterization of $W$ and $C(Y)$ in terms of the affine algebraic variety underlying $Z$.  For a simple example, we go back to $X^n_1$. For the Stenzel metric we know the tangent cone at infinity is $X^n_1$ itself. Since we are reversing the direction here, it does admit a weighted tangent cone at infinity isomorphic to $X^n_\infty$, so a priori $X^n_1$ could admit a Calabi-Yau metric with tangent cone at infinity given by $X^n_\infty$. This is also suggested by the construction of \cite{HN}.

\

\section{Appendix: Futaki and Matsushima theorem for polarized affine varieties}
In this section, we denote by $Z$ a tangent cone in the setting of Section 2.3.  The goal here is to prove Proposition \ref{prop3-5}, and some other related results. 
Recall we have proved that $Z$ is a polarized affine variety, endowed with a weak Ricci-flat K\"ahler cone metric. We fix an equivariant embedding of $Z$ into $\C^N$. The Reeb vector field $\xi_0$ generates a holomorphic action of a compact torus $\mathbb T$ on $Z$, which fixes the vertex of $Z$. Moreover, the action extends to  $\C^N$, through an embedding of $\mathbb T$ into the standard diagonal torus $\T^N$. In particular, if we denote by $\ct$ the Lie algebra of $\mathbb T$, then $\ct$ is naturally a subspace of $\bR^N$. Denote $\t^+=\t\cap (\bR^{+})^N$, then  $\xi_0\in \t^+$.

By Lemma 2.5 in \cite{HS}, there is a smooth family of $\T^N$-invariant K\"ahler cone metrics $\omega_\xi$ on $\C^N\setminus \{0\}$, parametrized by $\xi\in (\bR^+)^N$, such that for $\xi_1=(1, \cdots, 1)$, $\omega_{\xi_1}$ is the standard flat metric on $\C^N$, and for all  $\xi$, $\omega_\xi$ has Reeb vector field $\xi$ (called the type I deformation of $\omega_{\xi_1}$). Being a cone we have  $\omega_\xi=\frac{1}{4}dd^c r^2$, where $r$ is the distance function to the vertex with respect to $\omega_\xi$. For all $\xi$, the link $\{r=1\}$ is identified the unit sphere $S^{N-1}$ in $\C^N$, with the standard CR structure. 
For $\xi\in \t^+$, $\omega_\xi$ restricts to a $\mathbb T$ invariant K\"ahler cone metric on $Z$. Let $Y=Z\cap S^{N-1}$.  Notice $Y$ is in general different from, but naturally homeomorphic, to the link of $Z$ with respect to the Ricci-flat cone metric $\hat \omega$. 

We define 
$$V(\xi)=\int_Z e^{-r^2/2}(dd^c r^2)^n. $$
Up to multiplication by a dimensional constant, $V(\xi)$ is the same as the volume of $Y$ computed using the restriction of the metric $\omega_\xi$. For simplicity of notation we will denote the measure $d\mu=e^{-r^2/2}(dd^c r^2)^n$. Let $\eta=d^c \log r$ be the dual one-form of $\xi$. It is $\mathbb T$-invariant, and satisfies $\L_{r\p_r}\eta=0$.

\begin{lem} \label{lemA-1}
\begin{equation} \label{eqnA-1}
 dV(\delta\xi)=-n \int_Z \eta(\delta\xi)  d\mu
\end{equation}
\begin{equation}\label{eqnA-2}
\text{Hess}V(\delta\xi, \delta'\xi)=n(n+1)\int_Z \eta(\delta\xi)\eta(\delta'\xi) d\mu.
\end{equation}
In particular, 
$V(\xi)$ is strictly convex on $\t^+$. 
\end{lem}

This is proved in \cite{MSY} under the assumption that $Z\setminus\{0\}$ is smooth. We will perform the  calculation on the cone $Z$, from which it is evident that the appearance of singularities does not cause essential difficulties.

  We work on $\C^N\setminus\{0\}$, and denote the variation by $\delta (r^2)=r^2\phi$.  
Taking the variation of the equation $\L_{r\p_r}r^2=2r^2$, we obtain
\begin{equation}\label{eqnA-3}
d^c\phi(\xi)=-2\eta(\delta\xi).
\end{equation}
By definition the right hand side is radially invariant. It follows that $|\phi(r)|\leq C|\log r|$. Similarly  $|\nabla\phi(r)|\leq Cr^{-1}|\log r|$. 
We  compute the first variation
\begin{eqnarray*}
dV(\delta\xi)
&=& \int_Z e^{-r^2/2}(-\frac{1}{2}r^2\phi (dd^c r^2)^n+ndd^c(r^2\phi)(dd^c r^2)^{n-1})\\
&=&-\int_Z  \frac{1}{2}r^2\phi d\mu+\frac{n}{2} d(r^2)d^c(r^2\phi) e^{-r^2/2} (dd^c r^2)^{n-1}
\end{eqnarray*}

The second equality involves integration by parts. This can be verified by lifting to an $\mathbb T$-equivariant resolution $Z'$ (see \cite{Kollar} for the existence of such a resolution), and using the above estimate of $\phi$ and $|\nabla\phi|$ on $\C^N\setminus\{0\}$. 

Note for one-forms $\alpha$, $\beta$ on $\C^N\setminus\{0\}$, we have
\begin{equation} \label{eqnA-4}
n\alpha\wedge\beta\wedge (dd^c r^2)^{n-1}=\frac{1}{4}\langle \alpha, J\beta\rangle (dd^c r^2)^n. 
\end{equation}
Applying (\ref{eqnA-4}) with $\alpha=dr$, and $\beta=d^c(r^2\phi)$, we obtain
\begin{equation*}
dV(\delta\xi)=\frac{1}{4}\int_Z  d^c\phi(\delta\xi)r^2 d\mu.
\end{equation*}
Using (\ref{eqnA-3})  and the fact that $\eta(\delta\xi)$ is $r$-invariant,  this proves (\ref{eqnA-1}).

Now consider a new variation $\delta'\xi$ and accordingly $\delta' (r^2)=r^2\psi$. Then 
\begin{eqnarray*}
&&-\frac{1}{n}\text{Hess}V(\delta\xi, \delta'\xi)\\
&=& \int_Z \frac{1}{2}d^c\psi(\delta\xi)d\mu-\int_Z\frac{1}{2}r^2\psi \eta(\delta\xi)  d\mu+\int_Z\eta(\delta\xi) e^{-r^2/2} n dd^c (r^2\psi) (dd^c r^2)^{n-1}\\
&=&\I+\II+\III
\end{eqnarray*}

As above we use integration by parts to get 
\begin{eqnarray*}
&&\III\\
&=&-n\int_Z d(\eta(\delta\xi)) \psi d^c (r^2) e^{-\frac{r^2}{2}}(dd^c r^2)^{n-1}-n\int_Z r^2d(\eta(\delta\xi)) d^c\psi  e^{-\frac{r^2}{2}} (dd^c r^2)^{n-1}\\&&+ \frac{n}{2}\int_Z\eta(\delta\xi) e^{-r^2/2} d(r^2) r^2 d^c\psi (dd^c r^2)^{n-1}+2n\int_Zr^2\psi \eta(\delta\xi) drd^cr (dd^c r^2)^{n-1},
\end{eqnarray*}
 Applying (\ref{eqnA-4}) we see the first term in $\III$ vanishes since $\eta(\delta\xi)$ is $r$-invariant, the third term equals 
$-\frac{1}{2}\int_X \eta(\delta\xi)\eta(\delta'\xi)r^2 d\mu$, and the last term in $\III$ equals $\frac{1}{2}\int_Zr^2\psi \eta(\delta\xi) d\mu$ . For the second term in $\III$, we write
$$r^2 d(\eta(\delta\xi))=d(r^2\eta(\delta\xi)) -\eta(\delta\xi) d(r^2).$$
Notice that $\L_{\delta\xi}(r^2\eta)=0$, so 
\begin{equation} \label{eqnA-5}
d(r^2\eta(\delta\xi))=-\iota_{\delta\xi} d(r^2\eta)=-\frac{1}{2}\iota_{\delta\xi} dd^cr^2.
\end{equation}
So applying (\ref{eqnA-3}) and (\ref{eqnA-4}) we see
$$-n\int_Z r^2d(\eta(\delta\xi))d^c\psi  e^{-r^2/2}(dd^c r^2)^{n-1}=-\frac{1}{2}\int_Z d^c\psi(\delta\xi)d\mu -\int_Z\eta(\delta\xi)\eta(\delta'\xi)d\mu $$
Therefore 
$$
\III=-\frac{1}{2}\int_Z d^c\psi(\delta\xi)d\mu -\int_Z\eta(\delta\xi)\eta(\delta'\xi)(1+\frac{1}{2}r^2)d\mu+\frac{1}{2}\int_Zr^2\psi \eta(\delta\xi) d\mu
$$
Adding together $\I, \II, \III$, and using the fact that  $\eta(\delta\xi)$ and $\eta(\delta'\xi)$ are $r$-invariant, we get (\ref{eqnA-2}).

Finally, to see $V(\xi)$ is strictly convex, it suffices to show that if $\eta(\delta\xi)$ vanishes on $X$, then $\delta\xi=0$. This follows from (\ref{eqnA-5}). 

\begin{lem} \label{lemA-2}
$V(\xi)/V_n$ is a rational function with rational coefficients in the components of $\xi$, where $V_n$ is the volume of the round sphere $S^{2n-1}\subset \C^n$. 
\end{lem}

This is again proved in \cite{MSY} under the assumption that $Z\setminus \{0\}$ is smooth. For $\xi\in \mathcal R$, we define the index character
$$F(\xi, t)=\sum_{\alpha\in \Gamma^*} e^{-\langle\alpha, \xi\rangle t} \dim \H_\alpha, $$
where $\H_\alpha$ denotes the space of holomorphic functions on $Z$ with weight $\alpha$ under the  $\mathbb T$ action. 

It is shown in \cite{CS} that there is an asymptotic expansion (for $|t|\ll 1$)
$$F(\xi, t)=\frac{a_0(\xi)(n-1)!}{t^{n}}+\frac{a_1(\xi)(n-2)!}{t^{n-1}}+...$$
where $a_0(\xi)>0$ is a rational function in $\xi$ with rational coefficients. In particular $a_0(\xi)$ depends smoothly on $\xi$. We claim that  $a_0(\xi)=c_nV(\xi)$ where $c_n$ is a universal dimensional constant. Since both functions are continuous it suffices to prove this for a rational vector $\xi$, in which case one can  use the Riemann-Roch theorem for orbifolds (or more precisely, \emph{Deligne-Mumford stacks}) to obtain that $a_0(\xi)=\frac{1}{(n-1)!}\int_V  c_1(L)^{n-1}$, where $V$ is the quotient orbifold and $L$ is the descended ample orbi-line bundle. It follows from a similar calculation as in the smooth case that the latter can be computed using Chern-Weil theory, and we get $a_0(\xi)=c_nV(\xi)$. This proves the lemma. \\

\begin{lem} \label{lemA-3}
There are an integer $l$ and  a parallel section $s$ of $K_Z^l$, such that $(s\otimes \bar s)^{1/l}=\hat\omega^n$. 
\end{lem}

Suppose $Z$ is the Gromov-Hausdorff limit of a sequence $(X_i, L_i^{a_i}, a_i\omega_i, p_i)$ in $\K(n, \kappa, V)$ with $a_i\rightarrow\infty$. For simplicity of notation we only prove the case $\lambda=-1$, so that $L_i=K_{X_i}$. The proof of the other cases is similar. By the main results of \cite{DS}, we may find an integer $l$ and $C>0$, and holomorphic sections $s_i\in H^0(X_i, K_{X_i}^{l})$ with $|s_i(p_i)|=1$ and  $\int_{X_i}|s_i|^2\leq C$.  By the gradient estimate for holomorphic sections (Proposition 2.1 in \cite{DS}) there is a constant $D>0$ such that  $|\nabla_{\omega_i} s_i|_{L^\infty}\leq D$. Now after we rescale the manifold $X_i$ by a factor $a_i$, while fixing the Hermitian metric on $K_{X_i}$ (determined by the volume form of $\omega_i$) and the corresponding Chern connection,  we have $|\nabla_{a_i\omega_i}s_i|_{L^\infty}\leq a_i^{-1}D$. Then it follows by passing to a subsequence that $s_i$ converges locally uniformly to a section $s$ of $K_{Z^{reg}}^{l}$ which, over the regular part $Z^{reg}$, is parallel with respect to the Chern connection defined by the volume form of $\hat\omega$. Multiplying by a constant we may assume $(s\otimes \bar s)^{1/l}=\hat \omega^n$. Similar to the proof of Proposition 4.15 in \cite{DS} this implies that $Z$ has log terminal singularities so $K_Z^l$ is a well-defined line bundle, with a global section $s$ satisfying $(s\otimes \bar s)^{1/l}=\hat\omega^n$. 

\begin{rmk}
This lemma is the only place where we need to restrict our study to the smaller set $\K(n, \kappa, V)$ rather than $\K(n, \kappa)$. We expect the lemma to hold in greater generality and we leave this for future study.
\end{rmk}

As in \cite{HS}, we focus our attention on a hyperplane section in $\t^+$. Since $s$ is parallel on $Z$, we have $\L_{\xi_0}s=i a s$ for some $a\in \bR$. Since $(s\otimes \bar s)^{1/l}=\omega^n$ and  $\L_{r\p_r}\omega=2\omega$, it follows that $a=nl$. Similarly, for any $\xi\in \t$,  $\L_{\xi} s=i c(\xi) s$ for some linear function $c: \t\rightarrow \bR$.  Since $\mathbb T$ acts on the bundle $K_Z^l$, it is  not hard to see that $c$ has rational coefficients. Now we define
$$H=\{\xi\in \t^+| c(\xi)=nl\}. $$
By Lemma \ref{lemA-2} ,$V|_{H}$ is also a rational function with rational coefficients. Denote $\omega=\omega_{\xi_0}$, and $h=-\log |s|_{\omega}^{2/l}$. Then 
$$Ric(\omega)=-i\p\bp \log \omega^n=-i\p\bp h. $$
Consider a tangent vector $\delta\xi$ of $H$.  By (\ref{eqnA-5}) we have $\delta\xi=\frac{1}{2}J\nabla (r^2\eta(\delta\xi))$. So
$$\L_{J\delta\xi} h=-\Delta(r^2\eta(\delta\xi))+\frac{2}{l}c(\delta\xi), $$
and 
\begin{equation} \label{eqnA-6}
dV(\delta\xi)=\frac{1}{2}\int_Z r^2\eta(\delta\xi) d\mu=-\frac{1}{2}\int_Z \L_{J\delta\xi}h d\mu.
\end{equation}

\begin{prop} [Futaki Theorem] \label{propA-4}
$\xi_0$ is a critical point of $V|_{H}$. 
\end{prop}

Given this Proposition, it follows that $\xi_0$ is a critical point of a set of polynomial equations with rational coefficients. Now Lemma \ref{lemA-1} implies the Hessian of $V|_{H}$ is non-degenerate, so $\xi_0$ is indeed an isolated critical point on $H\otimes \C$. 
Then Proposition \ref{prop3-5} follows from an observation in \cite{MSY}. For completeness we provide a detailed argument here. 

\begin{lem}
Suppose \emph{$X=(x_1, \cdots, x_r)\in \C^r$} is an isolated zero of a system of polynomial equations with rational coefficients, then each $x_i$ is an algebraic number. 
\end{lem}

Suppose this fails,  without loss of generality we may assume $\{x_1, \cdots, x_t\}$ is a maximal algebraically independent subset of $\{x_1, \cdots, x_r\}$. Then for any $\{x_1', \cdots, x_t'\}$ such that $\{x_1, x_1', \cdots, x_t, x_t'\}$ is algebraically independent we can find an element $\tau\in Gal(\C/\Q)$ such that $\tau(x_i)=x_i'$. Here $Gal(\C/\Q)$ denotes the group of field automorphisms of $\C$ that fix elements in $\Q$. Clearly for any fixed $\delta>0$ we may assume $|x_i-x_i'|< \delta$ for  all $i\leq t$. Now  let $g_{t+1}(x)$ be the minimal polynomial of $x_{t+1}$ over $\Q(x_1, \cdots, x_{t})$. If we choose $x_i'(i=1, \cdots, t)$ as above, then we can find $x_{t+1}'\in\C$ such that  $|x_{t+1}'-x_{t+1}|=\epsilon(\delta)$ and $\tau(g_{t+1})(x_{t+1}')=0$, where $\epsilon(\delta)$ tends to zero as $\delta$ tends to zero. Now since $g_{t+1}(x_{t+1})=0$, we can choose $\sigma_{t+1}\in Gal(\C/\Q(x_1, \cdots, x_{t}))$ such that $\tau_{t+1}=\tau \circ \sigma_{t+1}$ sends $x_{t+1}$ to $x_{t+1}'$. Then we can proceed by induction to find for all $j\geq t+2$, an $x_j'\in \C$ with $|x_j'-x_j|\leq \epsilon(\delta)$, an element $\sigma_j\in Gal(\C/\Q(x_1, \cdots, x_{j-1}))$ such that $\tau_{j}=\tau_{j-1}\circ \sigma_j$ sends $x_j$ to $x_j'$.  It follows that $X'=(x_1', \cdots, x_r')$ is also a zero of the same system of polynomial equations. Let $\delta\rightarrow0$, we see that $X$ is not an isolated zero, contradiction.

\

If $Z$ is smooth, then the expression (\ref{eqnA-6})  is the usual Futaki-invariant adapted to K\"ahler cones. The crucial fact is that this is independent of the choice of the K\"ahler cone metric on $Z$ with fixed Reeb vector field. Hence we can compute it using the Ricci-flat cone metric, and derive the vanishing of (\ref{eqnA-6}). 
In general $Z$ is singular. We will use the results of pluripotential theory to prove Proposition \ref{propA-4}. 

\

  Notice that on $Y$ we have a Reeb foliation by $\xi_0$, a contact 1-form $\eta$, and  a transverse K\"ahler structure $\omega^T=\frac{1}{2}d\eta$ (strictly speaking, a K\"ahler current near the singular part of $Y$), all induced from $S^{N-1}$.  Let $\H$ be the space of bounded transverse K\"ahler potentials, i.e. the space of basic (i.e. $\mathbb T$-invariant), bounded, upper semi-continuous functions on $Y$ that is transversely pluri-subharmonic with respect to $\omega^{T}$. As usual such a K\"ahler potential $\phi$ gives rise to a transverse Monge-Amp\`ere measure, which together with the form $\eta$, defines a $\mathbb T$-invariant measure $(d\eta+dd^c\phi)^{n-1}\wedge \eta$ on $Y$. In the smooth case this agrees with the Riemannian volume form of the Sasaki structure defined by $\eta+d^c\phi$.

On the other hand, the above holomorphic section $s$ on $Z$ defines a volume form $\Omega$ on the smooth part of $Y$ by $(s\otimes \bar s)^{1/l}|_Y=dr\wedge \Omega$. So $\Omega$ determines a $\mathbb T$-invariant measure on $Y$, which we also denote by $\Omega$. An element $\phi$ in $\H$ then defines a $\mathbb T$-invariant measure $\Omega_\phi=e^{-\phi}\Omega$ on $Y$.

Let $\hat r$ be distance function to the vertex, defined by the metric $\hat \omega$. Write $\hat r=re^{\phi}$ for some $\mathbb T$-invariant function $\phi$ on $Z$, then the fact that  $\hat \omega$ and $\omega$ have the same Reeb vector field implies that $\phi$ is also $r$-invariant. So we may view $\phi$ as an element in $\H$. One then checks that   $(d\eta+dd^c\phi)^{n-1}\wedge \eta=C\Omega_{\phi}$ for a positive constant $C$. So it defines a \emph{weak transverse K\"ahler-Einstein metric}.

As in \cite{CDS3}, we define the \emph{Ding functional} 

$$\mathcal D(\phi)=I(\phi)-\log \int_Y\Omega_\phi$$

where 
$$I(\phi)=-\frac{1}{nV(\xi_0)}\sum_{i=0}^{n-1}\int_Y \phi (d\eta)^i\wedge (d\eta+dd^c\phi)^{n-1-i}\wedge \eta,$$
and the terms are made sense in terms of the usual pluri-potential theory.

 Now  given $\delta\xi$, let $f_t$ be the family of holomorphic transformations of $Z$ generated by $J\delta\xi$, and we denote by $\phi(t)$ the corresponding family of transverse K\"ahler potentials.  Then a direct calculation (similar to Lemma 12 in \cite{CDS3}) shows that

\begin{equation} \label{eqnA-7}
\frac{d}{dt}\mathcal D(\phi(t))=-\int_Z \L_{J\delta\xi}h d\mu 
\end{equation}
Notice the right hand side is independent of $t$. Just as in \cite{Bern}, given $\phi_0, \phi_1\in\H$,  one can find a bounded geodesic $\phi(t) (t\in [0, 1])$ in $\H$ connecting $\phi_0$ and $\phi_1$. The key property we need is

\begin{prop} \label{propA-5}
$\mathcal D$ is convex along $\phi(t)$. 
\end{prop}

It is straightforward to check that in our setting $\phi$ is a critical point of $\mathcal D$.  Then Proposition \ref{propA-5} implies that $\mathcal D$ is bounded below on $\H$. Then Proposition \ref{propA-4} follows from (\ref{eqnA-7}) and (\ref{eqnA-6}).

\

Therefore we are finally reduced to prove Proposition \ref{propA-5}. We also state two related results that is used in Section 3.  Let $\Aut(Z)$ be the group of holomorphic transformations of $Z$ that preserves $\xi_0$; in the notation of Section 3, this is a subgroup of $G_{\xi_0}$ that fixes $[Z]$ in $\Hilb$. 
The following results were proved in  \cite{Bern} and \cite{CDS3} for K\"ahler-Einstein $\Q$-Fano varieties.

\begin{prop} [Bando-Mabuchi theorem] \label{propA-7}
Ricci-flat K\"ahler cone metric on $Z$ with Reeb vector field $\xi_0$ is unique up to the action of the identity component of $\Aut(Z)$. 
\end{prop}

\begin{prop} [Matsushima theorem]\label{propA-6}
$\Aut(Z)$ is reductive. 
\end{prop}

\

Proposition \ref{propA-5}, \ref{propA-7} and \ref{propA-6} can be proved using  arguments analogous to the appendix of \cite{CDS3}, with the main technical  input from \cite{Bern} and \cite{BBEGZ}. We will only sketch below the key points that require extra care in our setting.   

\begin{rmk}
In the three dimensional case, we proved in \cite{DS} that $Y$ is a five dimensional Sasaki-Einstein orbifold, hence in that case all the above results can be alternatively obtained by direct computations similar to the case of smooth Sasaki-Einstein manifolds. 
\end{rmk}

 \
 
(1).  One difference in our setting is that we do not have a ``resolution of singularities" for Sasaki manifolds or affine cones. Notice Lemma \ref{lemA-3} implies that $Z$ has log terminal singularities. We can find a $\T^{\C, k}$-equivariant log resolution of singularities  $\pi: Z'\rightarrow Z$ (c.f. \cite{Kollar}), with simple normal crossing exceptional divisors $E_i$. So $K_{Z'}=\pi^*K_Z+\sum_i a_i E_i$ with $a_i>-1$ for all $i$. Let $E_1, \cdots, E_s$ be the set of  exceptional divisors that do not lie over the vertex (the other exceptional divisors are irrelevant). 
By construction we may assume the resolution is obtained by a sequence of blow-ups of the ambient space $\C^N$ at smooth $\T^\C$-invariant subvarieties. Let $P$ be the corresponding ambient space after blowing up and let $Y'=\pi^{-1}(Y)$. It follows that $Y'$ is naturally a smooth submanifold of $P$.  Let $\eta'=\pi^*\eta$, $\xi'=\pi^*\xi$, and $\omega'=\pi^*\omega$, then we obtain an induced foliation on $Y'$, and $(\eta', \xi', \omega')$ is a degenerate Sasaki structure on $Y'$.  It is in general not possible to deform this to a genuine Sasaki structure. 

But for our purpose we only need to deal with the transverse geometric properties of the foliation. By general theory, we find for all $i=1, \cdots, s$, a rational number $a_i>0$ and a Hermitian metric $h_i$ on the transverse holomorphic line bundle $E_i|_{Y'}$ with curvature form $\omega_i$, such that $\omega'_\epsilon=\omega'-\epsilon \sum_i a_i \omega_i$ is a transverse K\"ahler form on $Y'$ for  all $\epsilon>0$ sufficiently small. We fix such a $\epsilon\in \textbf{Q}$. 

We write $-K_{Z'}=-\pi^*K_Z-E+\Delta'$, where $E$ and $\Delta$ are both effective, $E$ has integer coefficients and $\Delta'$ has coefficients in $(0,1)$.  For simplicity we denote by $K_{Y'}$ and $K_Y$ the transverse canonical line bundles on $Y'$ and $Y$ respectively. Then we have $-K_{Y'}=-\pi^*K_Y-E+\Delta'$. Let $L=K_{Y'}^{-1}\otimes E$, then $L$ is isomorphic to $-\pi^*K_Y+\Delta'$. \\

\noindent (2).  We need a version of ``transverse Hodge decomposition theorem" for basic forms, i.e. forms $\alpha$ on $Y'$ satisfying $\iota_{\xi'}\alpha=0$ and $\mathcal L_{\xi'} \alpha=0$. One can define a transverse Hodge $*$ operator acting on basic forms, using the transverse volume form $\omega'_\epsilon$. Globally we use the $L^2$ inner product defined by $\omega'_\epsilon$ and $\eta'$.  Using the fact that $d\eta'$ is basic, one sees that the formal adjoint $d^*$ of $d$ is indeed given by $-*d*$. Then it is easy to work locally in the leaf space and develop the relevant elliptic theory for the basic Laplacian operator. One can also work out the analogue for $\bp$ operator. \\

\noindent (3).  One needs to check the local construction of pluripotential theory works well in our setting.  For example, we need  to approximate bounded pluri-subharmonic functions by a decreasing sequence of smooth functions which are almost pluri-subharmonic.  The results of \cite{Bern} use the construction of Blocki-Kolodziej, which depends on the choice of cut-off functions.   Notice we do not have $\mathbb T$-invariant cut-off functions on $Y'$ in general, but we can first do the construction using an arbitrary cut-off function, then take average over $\mathbb T$.\\

\noindent  (4). We need a version of the Kawamata-Viehweg vanishing that $H^{n,1}(Y', L)=0$. Note we may write $L=(-\pi^*K_Y-\epsilon \sum_i a_i E_i)+\Delta''$, such that $-\pi^*K_Y-\epsilon \sum_i a_i E_i$ admits a  Hermitian metric of  positive transverse curvature and $\Delta''$ still has coefficients in $(0, 1)$. Then we may apply the proof of Demailly \cite{Demailly}.  For the convenience of readers we provide here a detailed analytic proof in the case of compact K\"ahler manifolds, from which it is straightforward to extend to our setting, using the above transverse Hodge theory.

\begin{lem}
Let $X$ be a compact $n$ dimensional K\"ahler manifold and $L$ be a holomorphic line bundle over $X$. Suppose we can write $L^{\otimes k}=L'\otimes [F]$, where $L'$ is ample,  and $[F]$ is the line bundle defined by an effective divisor $F=\sum c_iF_i$ with normal crossing support and $k^{-1} c_i\in (0, 1)$. Then we have $H^{n, q}(X, L)=0$ for any $q\geq1$. 
\end{lem}

To prove this we choose a smooth Hermitian metric $h'$ on $L'$ with curvature $\omega>0$. Fix defining sections $s_i$ of $F_i$. These define a singular Hermitian metric $h_F$ on $[F]$ which is smooth away from $\cup F_i$, and with curvature $\sum_i a_i \delta_{F_i}$, where $\delta_{F_i}$ is the current of integration along $F_i$. Together with $h'$ this defines a singular Hermitian metric $h_0$ on $L$ with $iF_{h_0}\geq k^{-1}\omega$ as currents.  Given a smooth Hermitian metric $\tilde h$ on $[F]$,  for $\epsilon\in (0, 1]$  we obtain a smooth Hermitian metric 
$h_F (1+\epsilon h_F\tilde h^{-1})$ on $[F]$. Together with $h'$ this gives rise to a family of Hermitian metrics $h_\epsilon$ on $L$, that increase to $h_0$ as $\epsilon$ tends  to zero. Then a calculation (c.f. Lemma 16, \cite{CDS3}) shows that $iF_{h_\epsilon}\geq k^{-1}\omega-f_\epsilon^2\omega$ for a smooth function $f_\epsilon$ satisfying $0\leq f_\epsilon\leq C$ and $f_\epsilon$ converges to $0$ uniformly on any compact subset of $X\setminus \cup F_i$.  

Given $u\in \Omega^{n, q}(X, L)$  with $\bp u=0$, by the Kodaira-Nakano formula (\cite{Demailly2}) we have 
$$q^{-1}(\Delta_{\bp} u, u)_{\epsilon}\geq k^{-1}||u||_\epsilon^2-||f_\epsilon u||_{\epsilon}^2, $$
where the subscript $\epsilon$ denotes the $L^2$ inner product is defined in terms of $\omega$ and $h_\epsilon$. 
By standard elliptic theory, the operator $q^{-1}\Delta+f_\epsilon^2$ has an inverse $G_\epsilon$ with $||G_\epsilon u||_{\epsilon}^2\leq k||u||_{\epsilon}^2\leq kq^{-1}||u||_0^2$. So we can write 
$$u=q^{-1}\bp\bp^* G_\epsilon u+q^{-1}\bp^*\bp  G_\epsilon u+f_\epsilon^2 G_\epsilon u$$
Since $||G_\epsilon u||_{\epsilon}^2$ is uniformly bounded and $f_\epsilon\leq C$, it follows that $||\bp^*G_\epsilon u||_{\epsilon}^2$ is uniformly bounded. This implies that $||\bp^*G_\epsilon u||_{1}^2$ is also uniformly bounded.  By passing to a subsequence we may assume that as $\epsilon\rightarrow 0$, $\bp^*G_\epsilon u$ converges weakly to a limit $v$ in $L^2$. Since $\bp u=0$, we have
$$u=q^{-1}\bp\bp^*G_\epsilon u+\Pi(f_\epsilon ^2G_\epsilon u), $$
where $\Pi$ denotes the $L^2$ orthogonal projection to $Ker \bp$, defined in terms of the metric $h_\epsilon$.  Write $w_\epsilon=\Pi(f_\epsilon^2G_\epsilon u)$, then $||w_\epsilon||_\epsilon^2$ is uniformly bounded, so it converges  weakly to a limit $w$ in $L^2$. Moreover since $h_\epsilon$ converges to $h_0$ locally uniformly away from $\cup F_i$, $||w||_0^2<\liminf_{\epsilon\rightarrow0} ||w_\epsilon||_{\epsilon}^2$. It then follows that $u=q^{-1}\bp v+w$. We claim $||w||_0^2=\lim (w, w_\epsilon)_\epsilon$.  Indeed, writing $h_\epsilon=h_1 H_\epsilon$ $(\epsilon\in [0, 1])$ for a positive function $H_\epsilon$, then
$$(w, w_\epsilon)_\epsilon=\int \langle w, w_\epsilon\rangle_{h_1} H_\epsilon.$$
Since $||w_\epsilon||_{\epsilon}^2$ is uniformly bounded, we have $w_\epsilon H_{\epsilon}^{1/2}$ converges weakly in $L^2$ to $wH_0^{1/2}$. 
So
$$\int \langle w, w_\epsilon\rangle_{h_1} H_\epsilon^{1/2}H_0^{1/2}\rightarrow ||w||_0^2. $$
On the other hand, since $H_\epsilon\leq H_0$ and $H_\epsilon$ converges to $H_0$ away from $\cup F_i$, we have 
$$|\int  \langle w, w_\epsilon\rangle_{h_1} H_\epsilon^{1/2} (H_\epsilon^{1/2}-H_0^{1/2})|^2\leq ||w_\epsilon||_\epsilon^2\int |w|_{h_0}^2 H_0(H_\epsilon^{1/2}H_0^{-1/2}-1)^2\rightarrow 0.$$
This proves the claim.
Finally we have
$$(w, w_\epsilon)_\epsilon=(w, f_\epsilon^2 G_\epsilon u)_\epsilon\leq ||f_\epsilon^2 w||_{\epsilon}^2 ||G_\epsilon u||_\epsilon^2\leq ||f_\epsilon^2w||_0^2 ||G_\epsilon u||_\epsilon^2\rightarrow0,$$
  where the last inequality uses the fact that $f_\epsilon$ converges to zero uniformly on compact subset of  $X\setminus \cup F_i$. So $w=0$, and $u=q^{-1}\bp v$. 

\

Simons Center for Geometry and Physics, Stony Brook, U.S.A 

$\&$

Department of Mathematics, Imperial College London, U.K. 

Email: s.donaldson@imperial.ac.uk.\\

Department of Mathematics, Stony Brook University, U.S.A. 

Email: song.sun@stonybrook.edu\\

\end{document}